\documentclass{article}
\usepackage[left=1in,right=1in]{geometry}
\usepackage[utf8]{inputenc}
\usepackage{amsmath,amsfonts}
\usepackage{graphicx}
\usepackage[colorlinks=true, allcolors=blue]{hyperref}
\usepackage{subcaption}
\usepackage{comment} 
\usepackage{placeins}
\usepackage{booktabs}

\newcommand{\wt}{\widetilde}
\newcommand{\nvec}{n_{\mathrm{vec}}}

\title{Accurate Density of States Estimation: A Comparative Study}

\author{%
Dong Min Roh\thanks{Computational Research Division, Lawrence Berkeley National Laboratory, Berkeley, California 94720, United States}
\and
Chao Yang\footnotemark[1]
}

\date{Sep 2026}

\begin{document}

\maketitle

\begin{abstract}
We study density of states (DOS) estimation for large sparse real symmetric matrices by fitting the cumulative density of states (CDOS) and differentiating the fit. We use a stochastic Lanczos method to supply Ritz values and weights, which we average across Lanczos runs to construct CDOS midpoint data. Monotone piecewise-cubic interpolation of these data yields a smooth CDOS approximation that can be easily differentiated to yield a nonnegative, normalized, piecewise-quadratic DOS without requiring a Gaussian smoothing with a fixed  bandwidth in its construction. We also consider approximating CDOS by using Gaussian-process regression (GPR) with single- and double-Gaussian covariance kernels, yielding an approximation with uncertainty information. DOS approximation is then obtained by analytic derivative of the GPR-based CDOS. Experiments on matrices from diverse scientific applications are performed to compare these estimators with Gaussian-broadened stochastic Lanczos and Jackson-damped kernel polynomial approximations. At a common Gaussian validation resolution, we use error measures that assess local discrepancies, integrated errors, and agreement in the overall spectral shape. The results demonstrate that the midpoint spline can produce accurate approximations to the DOS with relatively few matrix--vector products, supporting the use of cumulative interpolation as a practical DOS estimation.

\end{abstract}

\section{Introduction}\label{sec:Intro}
Let $A\in\mathbb{R}^{n\times n}$ be symmetric, with eigenvalues
$\lambda_1\leq\lambda_2\leq\cdots\leq\lambda_n$ and corresponding orthonormal eigenvectors $u_1,\ldots,u_n$.
The spectral density of $A$ is the \emph{empirical spectral measure} defined by
\begin{equation}\label{eq:spectral_measure}
    \mu_A := \frac{1}{n}\sum_{j=1}^n \delta_{\lambda_j},
\end{equation}
where $\delta_{\lambda}$ denotes a unit point mass at $\lambda$.  It is also called the density of states (DOS) in physics and chemistry, where the eigenvalues of $A$ represent energy levels of a quantum system.

For every bounded function $f$ defined on the spectrum,
\begin{equation}\label{eq:spectral_measure_functional}
    \int_{\mathbb{R}} f(s)\,d\mu_A(s)
    = \frac{1}{n}\sum_{j=1}^n f(\lambda_j)
    = \frac{1}{n}\operatorname{tr} f(A),
\end{equation}
where $\operatorname{tr} (A)$ denotes the trace of the matrix $A$.
%This identity connects spectral-distribution estimation to stochastic trace estimation and polynomial approximation.

For any Borel set $B\subseteq\mathbb{R}$, $\mu_A(B)$ denotes the fraction
of eigenvalues of $A$ contained in $B$.  The \emph{cumulative density of
states} (CDOS) is the cumulative distribution function of $\mu_A$,
\begin{equation}\label{eq:cdos}
    c(t) := \mu_A(({-\infty},t])
    = \frac{1}{n}\sum_{j=1}^n H(t-\lambda_j),
\end{equation}
where
\begin{equation}\label{eq:Heavi}
    H(t) :=
    \begin{cases}
        1, & t\geq 0,\\
        0, & t<0.
    \end{cases}
\end{equation}
Thus $c(t)$ is nondecreasing, takes values in $[0,1]$, and equals the number of eigenvalues not exceeding $t$.

The DOS can also be defined as the distributional derivative of the CDOS, i.e.,
\begin{equation}\label{eq:dos}
    d(t) := c'(t)
    = \frac{1}{n}\sum_{j=1}^n \delta(t-\lambda_j),
\end{equation}
where $\delta(t-\lambda)$ is the Dirac delta distribution, characterized by
$\int f(t)\delta(t-\lambda)\,dt=f(\lambda)$ for smooth test functions $f$ \cite{lighthill1958introduction,dirac1981principles}.
Equations~\eqref{eq:spectral_measure}--\eqref{eq:dos} are three representations of the same spectral information: a measure, its cumulative distribution, and its distributional density.

The DOS is widely used in physics and chemistry to describe the global distribution of energy levels, while the CDOS directly answers eigenvalue-counting questions.
Spectral-distribution information is also useful in numerical linear algebra, including the design and analysis of iterative methods for linear systems and eigenvalue problems.
A complete evaluation of either quantity from its definition requires all eigenvalues of $A$, which is prohibitively expensive for a large matrix.
Scalable alternatives include the kernel polynomial method (KPM) and stochastic Lanczos quadrature \cite{kpmsurvey2006,LSY}.

\begin{figure}[tbp]
\centering
\includegraphics[width=\textwidth]{\detokenize{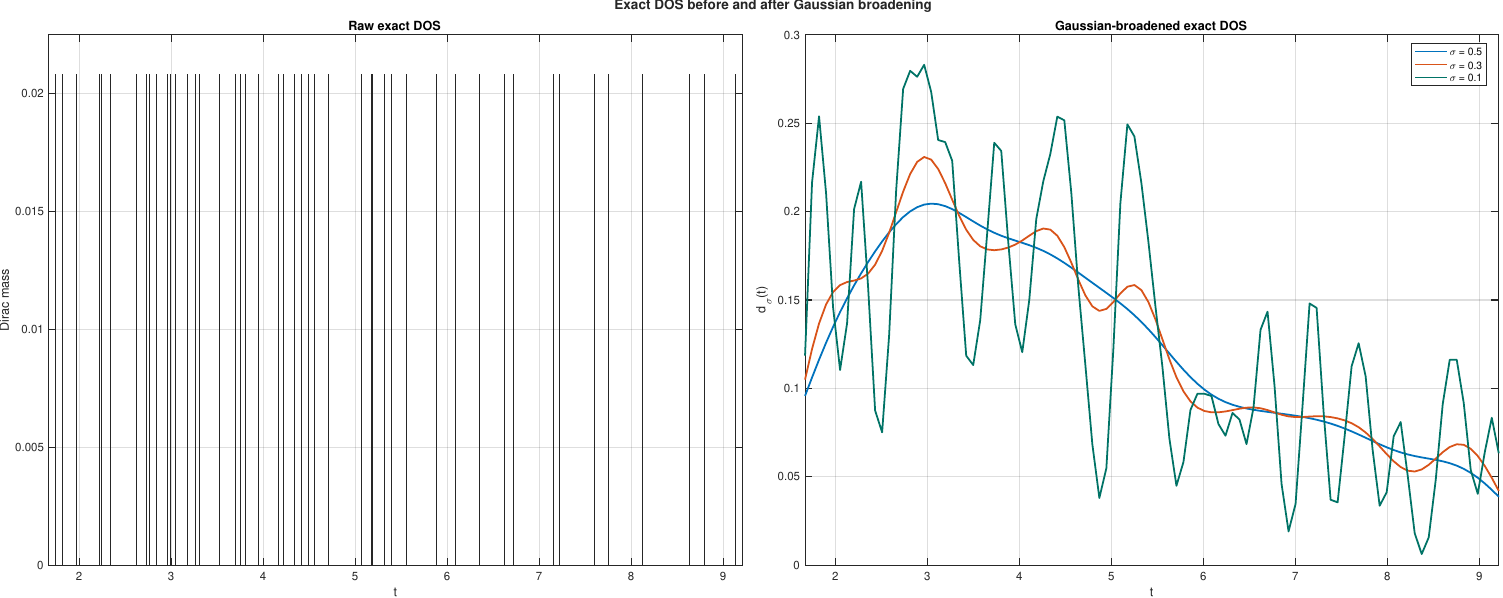}}
\caption{Exact and Gaussian-broadened DOS for the SuiteSparse
\texttt{mesh1e1} matrix ($n=48$)~\cite{davis2011university}. The left panel visualizes the normalized
exact DOS~\eqref{eq:dos} as impulses of mass $1/48$ at the eigenvalues.  The
right panel shows the Gaussian-regularized DOS~\eqref{eq:regDOS} for
$\sigma=0.5$, $0.3$, and $0.1$.  Decreasing $\sigma$ reveals progressively
finer spectral structure while making the curve more oscillatory.}
\label{fig:mesh1e1-DOS-regularization}
\end{figure}

Because the exact DOS is a distribution rather than an ordinary function, it is often regularized before numerical evaluation or visualization.
This can be achieved by convoluting it with a smooth, nonnegative kernel of unit integral, with Gaussian and Lorentzian kernels being common choices \cite{LSY}.
In this work, we focus on the Gaussian kernel
\begin{equation}\label{eq:gauss}
    g_{\sigma}(t)
    = \frac{1}{\sqrt{2\pi\sigma^2}}
      \exp\!\left(-\frac{t^2}{2\sigma^2}\right)
\end{equation}
with resolution parameter $\sigma>0$. Convolving the DOS with $g_\sigma$ yields the regularized DOS
\begin{equation}\label{eq:regDOS}
    d_{\sigma}(t)
    := (g_{\sigma}*d)(t)
    = \int_{\mathbb R}g_{\sigma}(t-s)\,d\mu_A(s)
    = \frac{1}{n}\sum_{j=1}^n g_{\sigma}(t-\lambda_j).
\end{equation}
In the distributional sense, $d_\sigma$ approaches the raw
DOS as $\sigma\to0$.
Figure~\ref{fig:mesh1e1-DOS-regularization} illustrates this regularization
using the $48\times48$ \texttt{mesh1e1} matrix from the SuiteSparse Matrix
Collection~\cite{davis2011university}. The impulse representation of the exact DOS identifies the
individual eigenvalue locations, but it is difficult to infer the overall
spectral distribution from this plot alone.  Gaussian broadening reveals
that distribution through the peaks of an ordinary function.  A larger
value of $\sigma$ produces a smoother curve with fewer, more easily
interpreted peaks, but suppresses fine spectral detail.  A smaller value
resolves more local structure but produces a more oscillatory, visually
noisy curve. Thus the selection of $\sigma$ is itself a
nontrivial choice of spectral resolution \cite{LSY}.

%The measure-based viewpoint suggests two routes to a smooth CDOS. One can first regularize a discrete Lanczos approximation to the DOS and integrate the resulting kernel expansion, or one can work directly with the cumulative distribution of the Lanczos quadrature measure and interpolate its staircase data. The second route has the advantage that monotonicity and the bounds $0\leq c(t)\leq1$ can be imposed explicitly. However, neither route is automatic: the Gaussian width affects the first estimator, generic interpolation may violate shape constraints, and Gaussian smoothing can introduce boundary bias when significant quadrature mass lies near a spectral endpoint.

Existing approaches for approximating the DOS include the kernel polynomial method (KPM) \cite{SilverRoder1994,kpmsurvey2006}, stochastic Lanczos quadrature \cite{LSY}, and maximum entropy methods \cite{DraboldSankey1993}. KPM uses polynomial expansions and stochastic trace estimates, whereas Lanczos quadrature represents spectral information through weighted Ritz values. An undamped polynomial truncation can produce negative values and Gibbs oscillations for KPM; Jackson damping addresses positivity through a positive reconstruction kernel, while introducing additional smoothing. A Gaussian-broadened Lanczos estimate is nonnegative, but may display pronounced quadrature-node peaks when the chosen bandwidth is too small relative to the available quadrature resolution.

Cumulative spectral approximation is also well established. Fischer develops interpolation of Lanczos-derived distribution data by monotone cubic functions \cite[Section~5.3]{fischer}; Lin, Saad, and Yang discuss spline-interpolated Lanczos CDOS in \cite[Appendix~C]{LSY}. Chen, Trogdon, and Ubaru analyze stochastic Lanczos approximation of the cumulative empirical spectral measure \cite{chen2021slq} and give a unified analysis of randomized Lanczos and polynomial quadrature methods \cite{chen2025matrixfree}. Building on these foundations, our focus is a DOS estimator based on averaged Lanczos midpoint data, followed by analytic differentiation, and an empirical assessment of the resulting accuracy--resolution tradeoff. The monotone interpolation formulas themselves are inherited from the established construction.

The exact finite-matrix CDOS is a staircase, not a smooth function. Nevertheless, its bounded, nondecreasing structure provides useful constraints for constructing a smooth surrogate, in contrast to the distribution-valued exact DOS. We compare monotone cubic interpolation with Gaussian-process (GP) regression using one or two covariance length scales. The spline enforces nonnegativity and unit mass of its derivative; the GP provides model-conditional uncertainty information but requires additional constraints to guarantee a valid density. Neither construction requires a Gaussian DOS-broadening bandwidth, unlike the stochastic Lanczos approach. We use a separate Gaussian bandwidth for the common-resolution validation of DOS in this paper.

The paper is outlined as follows.
Section~\ref{sec:CDOS_approx} develops the stochastic Lanczos framework and the two CDOS-based approximation methods. Section~\ref{sec:metrics} defines the validation resolution and error measures. Section~\ref{sec:results} presents the test matrices, comparison methods, experimental protocol, and matrix-by-matrix results, followed by a quantitative synthesis. Section~\ref{sec:conclusion} summarizes the findings and their limitations.

\section{Stochastic Lanczos and CDOS-Based DOS Estimation}\label{sec:CDOS_approx}

We first describe the spectral measure approximated by a single Lanczos run and then explain how random probes target the empirical measure. The resulting quadrature data support either direct Gaussian broadening or the midpoint-based cumulative fits developed below.

\subsection{Lanczos Quadrature Measure}\label{sec:lanczos_measure}

Let $v\in\mathbb{R}^n$ be a unit vector with eigenvector expansion
\begin{equation}\label{eq:v0_exp}
    v=\sum_{j=1}^n \alpha_j u_j,
    \qquad \sum_{j=1}^n\alpha_j^2=1.
\end{equation}
The vector $v$ induces the probability measure
\begin{equation}\label{eq:vector_measure}
    \mu_v:=\sum_{j=1}^n\alpha_j^2\delta_{\lambda_j}.
\end{equation}
Its action on a function and its cumulative distribution are
\begin{align}
    \int f(s)\,d\mu_v(s) &= v^T f(A)v,\label{eq:vector_measure_functional}\\
    F_v(t) &:= \mu_v(({-\infty},t])
    =\sum_{j=1}^n\alpha_j^2H(t-\lambda_j).\label{eq:spec_distr}
\end{align}
The associated distributional density is
\begin{equation}\label{eq:DOS_from_distr}
    d_v(t)=F_v'(t)
    =\sum_{j=1}^n\alpha_j^2\delta(t-\lambda_j).
\end{equation}
In general, $F_v$ is a vector-weighted spectral distribution rather than the CDOS $c$. 
Only when the weights are uniform, i.e., $\alpha_j^2=1/n$, the equality $F_v=c$ holds.
% Uniform weights $\alpha_j^2=1/n$ guarantee $F_v=c$. More generally, equality requires the weight in each eigenspace to equal its multiplicity divided by $n$.
% More generally, when eigenvalues are repeated, equality requires that the total weight in each eigenspace equal its multiplicity divided by $n$.

In exact arithmetic and before breakdown, an $M$-step Lanczos process starting from $q_1=v$ constructs an orthonormal basis
$Q_M=[q_1,\ldots,q_M]$ for the Krylov subspace
\begin{equation}\label{eq:Krylov}
    \mathcal{K}_M(A;v)
    =\operatorname{span}\{v,Av,\ldots,A^{M-1}v\}
\end{equation}
and a symmetric tridiagonal matrix $T_M$ satisfying
\begin{equation}\label{eq:lanfact}
    AQ_M=Q_MT_M+\beta_{M+1}q_{M+1}e_M^T,
    \qquad Q_M^TQ_M=I_M.
\end{equation}
The Lanczos basis vectors have the form
\begin{equation}\label{eq:basis_poly_expr}
    q_{j+1}=p_j(A)v,
    \qquad j=0,\ldots,M-1,
\end{equation}
where the $p_j$ are orthonormal with respect to
\begin{equation}\label{eq:inn_prod_GAL}
    \langle p,q\rangle_{\mu_v}
    :=\int p(s)q(s)\,d\mu_v(s)
    =v^Tp(A)q(A)v.
\end{equation}
Thus $T_M$ is the Jacobi matrix associated with the orthogonal polynomials for $\mu_v$ \cite{fischer,golubwelsch}.

Let
\begin{equation}\label{eq:ritz_pairs}
    T_My_k=\theta_k y_k,
    \qquad \|y_k\|_2=1,
    \qquad \tau_k=e_1^Ty_k,
    \quad k=1,\ldots,M.
\end{equation}
The Ritz values $\theta_k$ and positive weights $\omega_k:=\tau_k^2$ define the Gaussian quadrature measure
\begin{equation}\label{eq:quadrature_measure}
    \mu_{v,M}:=\sum_{k=1}^M\omega_k\delta_{\theta_k},
    \qquad \omega_k>0,
    \qquad \sum_{k=1}^M\omega_k=1.
\end{equation}
The Gaussian quadrature rule gives
\begin{equation}\label{eq:gauss_quadrature_exactness}
    \int p(s)\,d\mu_v(s)
    =v^Tp(A)v
    =e_1^Tp(T_M)e_1
    =\sum_{k=1}^M\omega_kp(\theta_k)
    =\int p(s)\,d\mu_{v,M}(s),
    \qquad p\in\Pi_{2M-1},
\end{equation}
where $\Pi_r$ is the space of polynomials of degree at most $r$ \cite{fischer}.
Thus $\mu_v$ and $\mu_{v,M}$ define the same integration functional on
$\Pi_{2M-1}$, or equivalently, have the same moments through degree
$2M-1$.  The two measures are not generally identical; they are
indistinguishable only by polynomial test functions in this space.

The DOS and CDOS representations now follow directly from the same quadrature measure.
Its discrete DOS is
\begin{equation}\label{eq:quadrature_dos}
    d_{v,M}(t)=\sum_{k=1}^M\omega_k\delta(t-\theta_k),
\end{equation}
whereas its discrete cumulative distribution is
\begin{equation}\label{eq:new_distr}
    F_{v,M}(t)=\sum_{k=1}^M\omega_kH(t-\theta_k).
\end{equation}
In the sense of distributions, $F_{v,M}'=d_{v,M}$.
The polynomial exactness in \eqref{eq:gauss_quadrature_exactness} does not apply directly to the discontinuous function $H(t-s)$; it supplies moment matching between the two measures.
For a nonzero cumulative error, moment matching forces alternating cancellation, as expressed by the classical sign-change result in \cite[Theorem~2.2.5]{fischer}. The interlocking behavior motivates, but does not by itself guarantee, accurate pointwise interpolation of the CDOS.

\subsection{Stochastic Approximation of the Empirical Spectral Measure}

There are at least two issues in using \eqref{eq:quadrature_dos} to approximate \eqref{eq:dos} directly. First, \eqref{eq:quadrature_dos} contains fewer $\delta$ peaks than the true DOS, so it cannot resolve all spectral features, even though the weights $\omega_k$ indicate the concentration of spectral mass near the nodes $\theta_k$. Second, the weights depend on the initial vector $v$ of the Lanczos process. For a general $v$, the coefficients $\alpha_j^2$ in \eqref{eq:vector_measure} are not uniform. Consequently, even a complete $n$-step Lanczos process recovers the vector-induced measure $\mu_v$, not necessarily the empirical spectral measure $\mu_A$. 
% Increasing $M$ therefore removes the quadrature error relative to $\mu_v$, but it does not remove the dependence on a single realized starting vector.
%For a general starting vector, the weights $\alpha_j^2$ in \eqref{eq:vector_measure} are not uniform.

Gaussian smoothing addresses the first issue, but introduces a resolution parameter; we return to it in Section~\ref{sec:CDOS}. To address dependence on the starting vector, we average the quadrature measures from independent random probes. Let
$g^{(1)},\ldots,g^{(\nvec)}$ be independent standard Gaussian vectors and define the normalized Lanczos starting directions by
\begin{equation}
    v^{(l)}=\frac{g^{(l)}}{\lVert g^{(l)}\rVert_2}.
\end{equation}
This normalization is required by the convention $q_1=v^{(l)}$ used above, because $q_1$ is the first column of an orthonormal Lanczos basis. It changes only the scale of the starting vector, not the Krylov subspace that it generates.
The vectors $v^{(l)}$ are unit vectors drawn isotropically from the unit sphere, and hence
\begin{equation}\label{eq:isotropic_probe}
    \mathbb{E}\bigl[v^{(l)}(v^{(l)})^T\bigr]=\frac{1}{n}I.
\end{equation}
Then, for every suitable $f$,
\begin{equation}\label{eq:stochastic_trace}
    \mathbb{E}\bigl[(v^{(l)})^Tf(A)v^{(l)}\bigr]
    =\frac{1}{n}\operatorname{tr}f(A)
    =\int f(s)\,d\mu_A(s).
\end{equation}
Hence $\mathbb{E}[\mu_{v^{(l)}}]=\mu_A$ in the weak sense.
This is the normalized-direction version of stochastic trace estimation. In the more familiar formulation with an unnormalized standard Gaussian vector,
\begin{equation}
    \mathbb{E}\bigl[(g^{(l)})^Tf(A)g^{(l)}\bigr]
    =\operatorname{tr}f(A).
\end{equation}
Thus division by $n$ targets the normalized trace associated with $\mu_A$. If the Lanczos process is run with the unit direction $v^{(l)}$, the corresponding raw-Gaussian estimator would multiply its quadrature rule by $\lVert g^{(l)}\rVert_2^2/n$. Here we instead use normalized directions directly, so that every sample measure has total mass one and its expectation is $\mu_A$ \cite{LSY,AvronToledo2011}.

We run $M$ Lanczos steps for each probe $v^{(l)}$, producing nodes $\theta_k^{(l)}$ and weights
$\omega_k^{(l)}=(\tau_k^{(l)})^2$.
The averaged quadrature measure is
\begin{equation}\label{eq:DOS_from_distr_avg}
    \widehat\mu_M
    :=\frac{1}{\nvec}\sum_{l=1}^{\nvec}
      \sum_{k=1}^M\omega_k^{(l)}\delta_{\theta_k^{(l)}}.
\end{equation}
Its distributional DOS and CDOS representations are, respectively,
\begin{align}
    \widehat d_M(t)
    &=\frac{1}{\nvec}\sum_{l=1}^{\nvec}
      \sum_{k=1}^M\omega_k^{(l)}\delta(t-\theta_k^{(l)}),
      \label{eq:stochastic_dos}\\
    \widehat c_M(t)
    &=\frac{1}{\nvec}\sum_{l=1}^{\nvec}
      \sum_{k=1}^M\omega_k^{(l)}H(t-\theta_k^{(l)}).
      \label{eq:stochastic_cdos}
\end{align}
The DOS measure has total mass one and is nonnegative; its CDOS takes values in $[0,1]$ and is nondecreasing.

To distinguish Lanczos quadrature error from stochastic sampling error, consider first completing the Lanczos process in exact arithmetic. Recall that $\mu_{v,M}$ in~\eqref{eq:quadrature_measure} denotes the quadrature measure obtained from $M$ Lanczos steps starting with $v$. Thus $\mu_{v^{(l)},n}$ denotes this measure for the $l$th random probe with $M=n$, where $n$ is the matrix dimension. At completion, the quadrature measure equals the exact vector-induced spectral measure:
\begin{equation}
    \mu_{v^{(l)},n}=\mu_{v^{(l)}}
    =\sum_{j=1}^n |u_j^Tv^{(l)}|^2\delta_{\lambda_j}.
\end{equation}
If Lanczos terminates exactly before $n$ steps, the completed rule already has this property, and we use $\mu_{v^{(l)},n}$ as shorthand for that rule. Completion therefore removes the quadrature error, but the weights still depend on the random probe and need not equal the uniform weights $1/n$ of $\mu_A$. Since $\mathbb{E}[|u_j^Tv^{(l)}|^2]=1/n$, averaging increasingly many independent probes recovers those uniform weights. Consequently, for each bounded continuous test function $f$, the strong law of large numbers gives
\begin{equation}\label{eq:stochastic_dos_limit}
    \int f(s)\,d\widehat\mu_n(s)
    =\frac{1}{\nvec}\sum_{l=1}^{\nvec}
       (v^{(l)})^Tf(A)v^{(l)}
    \xrightarrow[\nvec\to\infty]{\mathrm{a.s.}}
    \frac{1}{n}\operatorname{tr}f(A)
    =\int f(s)\,d\mu_A(s).
\end{equation}
Here ``a.s.'' means \emph{almost surely}: the convergence holds with probability one over the random probes as their number $\nvec$ tends to infinity, with the matrix dimension $n$ fixed. Equation~\eqref{eq:stochastic_dos_limit} states that $\widehat\mu_n$ converges weakly to $\mu_A$, or equivalently that $\widehat d_n$ converges to $d$ in the distributional sense. It is not a pointwise limit of ordinary DOS functions. This is a theoretical consistency statement, not a requirement to perform $n$ Lanczos steps in practice. 
With $M\ll n$ steps per probe, we observe that the approximation is generally good enough using only a few $\nvec$ random starting directions.
% With $M\ll n$ steps per probe, the approximation generally contains both Lanczos quadrature error and finite-sample error from using only $\nvec$ random starting directions.

%Finally, convolving \eqref{eq:stochastic_dos} with the Gaussian kernel gives the regularized Lanczos DOS
%\begin{equation}\label{eq:DOS_from_distr_reg}
%    \widehat d_{\sigma,M}(t)
%    =\frac{1}{\nvec}\sum_{l=1}^{\nvec}
%      \sum_{k=1}^M\omega_k^{(l)}g_{\sigma}(t-\theta_k^{(l)}).
%\end{equation}
%The next section develops three smooth CDOS approximations from these common Lanczos quadrature data.

\subsection{Smoothing the DOS Approximations}\label{sec:CDOS}
\begin{figure}[tbp]
\centering
\includegraphics[width=\textwidth]{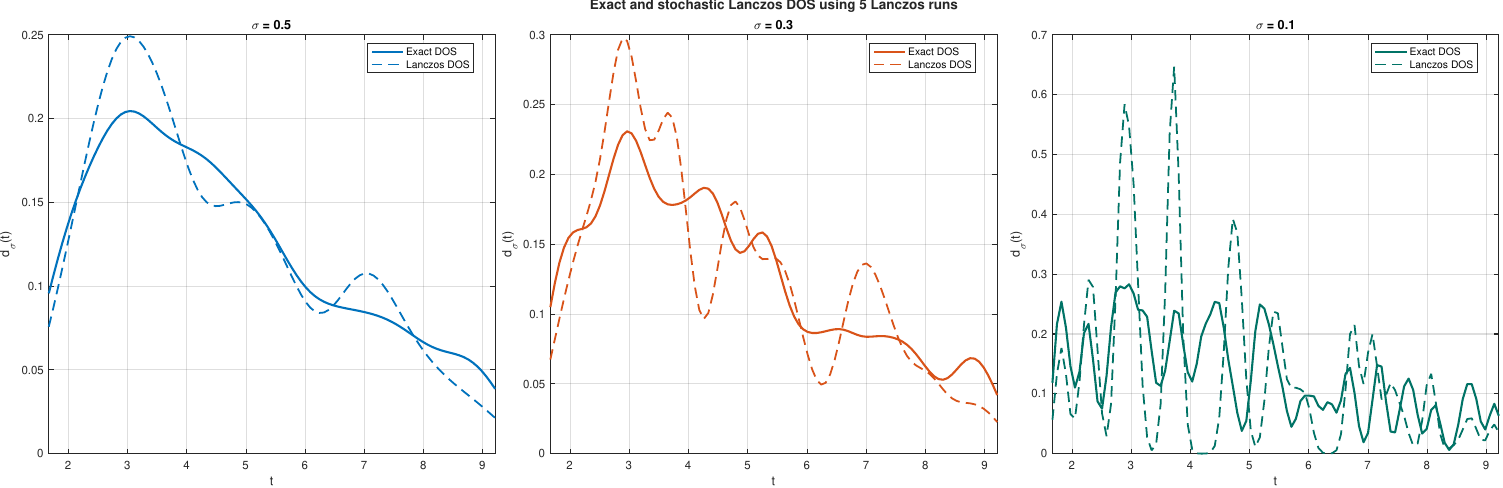}
\caption{Exact and stochastic Lanczos Gaussian-regularized DOS curves for
the SuiteSparse \texttt{mesh1e1} matrix~\cite{davis2011university}, using five independent Lanczos runs
with $M=10$ steps.  The exact DOS $d_\sigma$ is shown by a solid line and the
Lanczos approximation $\widehat d_{\sigma,M}$ by a dashed line.  From left
to right, the common smoothing bandwidth is $\sigma=0.5$, $0.3$, and $0.1$.
As the bandwidth decreases, the Lanczos approximation becomes substantially
more oscillatory than the exact regularized DOS.}
\label{fig:mesh1e1-DOS-Lanczos-sigma}
\end{figure}

%The stochastic Lanczos measure $\widehat\mu_M$ in \eqref{eq:DOS_from_distr_avg} has two exact discrete representations: the DOS $\widehat d_M$ and the staircase CDOS $\widehat c_M$.
One way to address the first issue raised in the previous section is to smooth the DOS approximation $\widehat d_M$ by convolving it with a Gaussian kernel, as in \eqref{eq:regDOS}. This yields the regularized Lanczos DOS
\begin{equation}\label{eq:DOS_from_distr_reg}
    \widehat d_{\sigma,M}(t)
    =\frac{1}{\nvec}\sum_{l=1}^{\nvec}
      \sum_{k=1}^M\omega_k^{(l)}g_{\sigma}(t-\theta_k^{(l)}).
\end{equation}
Choosing $\sigma$ remains difficult.  Figure~\ref{fig:mesh1e1-DOS-Lanczos-sigma}
compares the exact regularized DOS $d_\sigma$ with the stochastic Lanczos
DOS $\widehat d_{\sigma,M}$ using the same bandwidth for both curves.  For
$\sigma=0.5$, the broad Gaussian kernels suppress most of the differences
between the two underlying discrete measures, and the Lanczos approximation somewhat
captures the overall spectral shape.  At $\sigma=0.3$, discrepancies in the
locations and heights of individual peaks become more visible.  At
$\sigma=0.1$, both curves contain finer structure, but the Lanczos estimate
is substantially more oscillatory and develops sharp peaks and valleys that
are absent or much less pronounced in the exact curve.

This behavior follows from the different discrete measures being smoothed.
As $\sigma\to0$, $d_\sigma$ approaches the exact spectral measure supported
at all eigenvalues, whereas, for fixed $M$ and $\nvec$, the Lanczos curve
approaches the finite quadrature measure $\widehat\mu_M$.  Each Lanczos run
contributes only $M$ Ritz values, and when $\sigma$ is small relative to the
effective gaps between the quadrature nodes, neighboring Gaussian kernels
overlap too little to form a smooth curve.  Increasing $M$ supplies more
Ritz values and can improve the placement of spectral mass; increasing
$\nvec$ also reduces stochastic variation.  Nevertheless, direct Gaussian
smoothing still requires a bandwidth that is large enough for the available
quadrature resolution, potentially suppressing genuine fine-scale spectral
features.  This limitation motivates constructing a smooth CDOS first and
then differentiating it.

%This section constructs smooth approximations from the same Ritz values and quadrature weights.
%The first method convolves the measure with a Gaussian kernel and integrates analytically.
%The other two methods interpolate representative midpoint data extracted from a cumulative quadrature distribution.

\subsection{Midpoint Data for CDOS Approximation}
\label{sec:midpoint_data}

Instead of smoothing the DOS directly, we construct a smooth surrogate of
the nondecreasing CDOS. We
then differentiate the smooth CDOS approximation to obtain a DOS
approximation.  For an individual cumulative quadrature measure
$F_{v,M}$ in \eqref{eq:new_distr}, define the midpoint of its $j$th jump by
\begin{equation}\label{eq:interpol_pt}
    \vartheta_j
    :=\sum_{k=1}^{j-1}\omega_k+\frac{\omega_j}{2},
    \qquad j=1,\ldots,M.
\end{equation}
The pairs $(\theta_j,\vartheta_j)$ provide representative samples of the
unknown vector-weighted distribution $F_v$; this use of jump midpoints is
motivated by the interlocking properties of Gaussian quadrature
\cite[Section~5.3]{fischer}.  For the Chebyshev distribution, the exact
cumulative distribution passes through every jump midpoint, so that
$F_v(\theta_j)=\vartheta_j$ \cite[Example~2.2.3]{fischer}.  For a general
spectral measure, however, the sign-change theorem does not guarantee this
pointwise identity.

\begin{figure}[htbp]
\centering
\includegraphics[width=\textwidth]{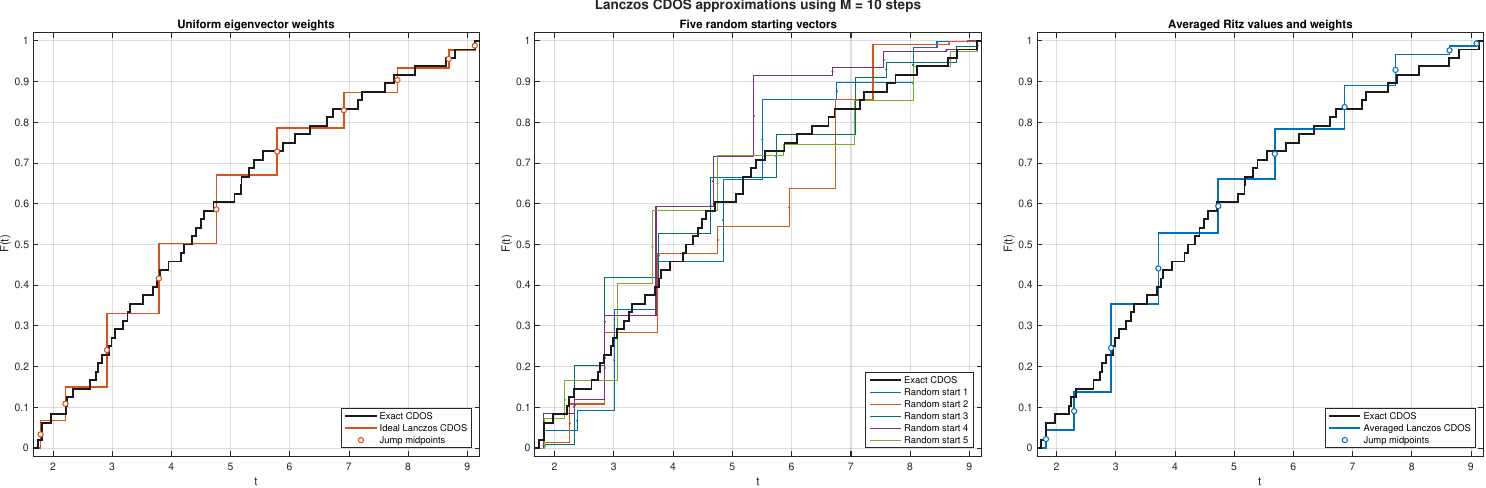}
\caption{Cumulative Lanczos quadrature measures for the SuiteSparse
\texttt{mesh1e1} matrix~\cite{davis2011university} using $M=10$ Lanczos steps. Left: the ideal
starting vector~\eqref{eq:ideal_starting_vector} gives uniform eigenvector
weights, and the jump midpoints of $F_{v_\star,M}$ closely track the exact
CDOS in this example.  Middle: five independent random starting directions
produce visibly different vector-weighted cumulative distributions whose
midpoints do not show the same close agreement.  Right: averaging corresponding ordered
Ritz values and their quadrature weights across those five runs produces
$\bar F_M$, whose midpoint data $(\bar\theta_j,\bar\vartheta_j)$ closely
follow the exact CDOS.}
\label{fig:mesh1e1-CDOS-Lanczos}
\end{figure}

Figure~\ref{fig:mesh1e1-CDOS-Lanczos} illustrates the construction of the interpolation data points that we use in three
steps, using the \texttt{mesh1e1} matrix and $M=10$ Lanczos steps.

\paragraph{Ideal starting vector.}
First consider the informative but impractical starting vector
\begin{equation}\label{eq:ideal_starting_vector}
    v_\star:=\frac{1}{\sqrt n}\sum_{j=1}^n u_j.
\end{equation}
Its eigenvector coefficients satisfy $\alpha_j^2=1/n$, and hence
$\mu_{v_\star}=\mu_A$ and $F_{v_\star}=c$.  The left panel of
Figure~\ref{fig:mesh1e1-CDOS-Lanczos} shows that the jump midpoints of the
corresponding Lanczos CDOS $F_{v_\star,M}$ closely track the exact CDOS in this example,
although they do not coincide exactly with it at every point.  We therefore call $v_\star$ the \emph{ideal}
starting vector.  Constructing it requires all eigenvectors of $A$, so it
is unavailable in a practical large-scale computation.

\paragraph{Individual random starting vectors.}
In practice, each Lanczos run begins with an independent random normalized
direction.  Its coefficients $\alpha_j^2$ are generally nonuniform, so the
resulting $F_{v,M}$ approximates the vector-weighted distribution $F_v$
rather than the empirical CDOS $c$.  The middle panel of
Figure~\ref{fig:mesh1e1-CDOS-Lanczos} shows five such runs.  Their staircase
distributions vary visibly, and their jump midpoints do not possess the
alignment seen in the ideal case.

\paragraph{Averaging across Lanczos runs.}
To reduce this probe-to-probe variation, we sort the Ritz values within each
run in increasing order, retaining the association between each node and its quadrature weight. For each position $j$, we then average the Ritz values and their associated weights across runs:
\begin{equation}\label{eq:new_interpol_pts}
    \bar{\theta}_j
    :=\frac{1}{\nvec}\sum_{l=1}^{\nvec}\theta_j^{(l)},
    \qquad
    \bar{\omega}_j
    :=\frac{1}{\nvec}\sum_{l=1}^{\nvec}\omega_j^{(l)},
    \qquad j=1,\ldots,M.
\end{equation}
This averaging procedure
is a compression of the full averaged measure $\widehat\mu_M$ in
\eqref{eq:DOS_from_distr_avg}; it is not algebraically identical to averaging
the measures, which generally produces up to $M\nvec$ distinct point masses.
Consequently, neither the polynomial exactness of the individual quadrature rules nor the unbiasedness of the uncompressed probe measures automatically carries over to the compressed measure. The construction here assumes that every run supplies the same number $M$ of strictly ordered nodes.

The compressed cumulative distribution is
\begin{equation}\label{eq:new_distr2}
    \bar F_M(t):=\sum_{j=1}^{M}\bar{\omega}_jH(t-\bar{\theta}_j),
\end{equation}
and the midpoint of its $j$th jump is
\begin{equation}\label{eq:avg_interpol_pt}
    \bar{\vartheta}_j
    :=\sum_{k=1}^{j-1}\bar{\omega}_k+\frac{\bar{\omega}_j}{2},
    \qquad j=1,\ldots,M.
\end{equation}

The right panel of Figure~\ref{fig:mesh1e1-CDOS-Lanczos} uses the averaged Ritz values
and weights from the same five random runs shown in the middle panel.  The
midpoints $(\bar\theta_j,\bar\vartheta_j)$ of the resulting compressed CDOS
closely follow the exact CDOS and thereby mimic the useful behavior of the
ideal case.  This alignment is an empirical observation rather than a
general identity guaranteed by Gaussian quadrature; it motivates using the
averaged midpoint pairs as training data for the smooth CDOS approximations
below.

\subsection{CDOS and DOS Approximation by Monotone Cubic Interpolation}
\label{sec:spline}

We construct a piecewise-cubic interpolant through the averaged midpoint
data $(\bar\theta_j,\bar\vartheta_j)$, assuming $M\geq2$. To impose the limiting CDOS values,
we augment these data with artificial endpoints
$(\bar\theta_0,\bar\vartheta_0)$ and $(\bar\theta_{M+1},\bar\vartheta_{M+1})$, following \cite{fischer}.
We set
\begin{equation}\label{eq:bnd_pt1}
    \bar{\theta}_0 := 
    \begin{cases}
        \bar{\theta}_1/2, & \bar{\theta}_1>0,\\
        \bar{\theta}_1-(\bar{\theta}_2-\bar{\theta}_1)/2,
        & \bar{\theta}_1\leq0,
    \end{cases}
\end{equation}
and
\begin{equation}\label{eq:bnd_pt2}
    \bar{\theta}_{M+1} := 
    \begin{cases}
        \bar{\theta}_M/2, & \bar{\theta}_M<0,\\
        \bar{\theta}_M+(\bar{\theta}_M-\bar{\theta}_{M-1})/2,
        & \bar{\theta}_M\geq0.
    \end{cases}
\end{equation}

We set $\bar\vartheta_0=0$ and $\bar\vartheta_{M+1}=1$. These artificial endpoints close the interpolant to unit mass; they are heuristic extrapolations, not certified bounds on the spectrum. In particular, their sign-dependent definitions depend on the origin of the spectral variable and can affect the estimate near the boundaries.

We connect consecutive data points with a monotone piecewise-cubic Hermite interpolant, denoted by $\wt{c}_{cs}(t)$. We use the term ``spline'' for this $C^1$ construction; it is not in general a $C^2$ cubic spline.
Following \cite{fischer, fritsch1980monotone, fritsch1984method}, we obtain for $t\in[\bar{\theta}_j, \bar{\theta}_{j+1}]$
\begin{equation}\label{eq:CDOS_CS}
\begin{split}
    \wt{c}_{cs}(t) = & \left(\frac{d_j+d_{j+1}-2\Delta_j}{h_j^2}\right)(t-\bar{\theta}_j)^3 \\
    & +\left(\frac{-2d_j-d_{j+1}+3\Delta_j}{h_j}\right)(t-\bar{\theta}_j)^2 + d_j(t-\bar{\theta}_j) + \bar{\vartheta}_j,
\end{split}
\end{equation}
where $h_j:=\bar{\theta}_{j+1}-\bar{\theta}_j$ and $\Delta_j:=(\bar{\vartheta}_{j+1}-\bar{\vartheta}_j)/h_j$ for $j=0,1,\ldots,M$.
We follow \cite{fischer, brodlie1980review, fritsch1984method} for the choice of the $d_j$'s
\begin{equation}\label{eq:djs}
    d_j:=
    \begin{cases}
        \frac{\Delta_{j-1}\Delta_j}{\xi_j\Delta_j+(1-\xi_j)\Delta_{j-1}} &\mbox{ for } \Delta_{j-1}\Delta_j>0, \\
        0 &\mbox{ otherwise, }
    \end{cases}
    \quad j=1,2,\ldots,M,
\end{equation}
where $\xi_j:=(h_{j-1}+2h_j)/(3(h_{j-1}+h_j))$.
The boundary conditions $d_0$ and $d_{M+1}$ are chosen to be
\begin{equation}\label{eq:d0dM+1}
    \begin{cases}
        d_0 &= \frac{h_0}{h_1}(3\Delta_1-(2d_1+d_2))+3\Delta_0-2d_1, \\
        d_{M+1} &= \frac{h_M}{h_{M-1}}(3\Delta_{M-1}-(2d_M+d_{M-1})) + 3\Delta_M - 2d_M.
    \end{cases}
\end{equation}
% Author check before submission: the uploaded mpci.m uses delta(n-1),
% rather than delta(n-2), in the first term of its right-endpoint slope.
% This differs from the Delta_{M-1} term above and from Fischer (5.3.4).
% Confirm which implementation generated the reported figures before
% changing the formula or recomputing any results.
To preserve monotonicity of $\wt{c}_{cs}(t)$ in the subintervals $[\bar{\theta}_0, \bar{\theta}_1]$ and $[\bar{\theta}_M, \bar{\theta}_{M+1}]$, we set $d_0=0$ or $d_0=3\Delta_0$ if $d_0$ in \eqref{eq:d0dM+1} is negative or bigger than $3\Delta_0$, respectively.
Similarly, we set $d_{M+1}=0$ or $d_{M+1}=3\Delta_M$ if $d_{M+1}$ in \eqref{eq:d0dM+1} is negative or bigger than $3\Delta_M$, respectively. 
On $[\bar\theta_0,\bar\theta_{M+1}]$, the resulting function is continuously differentiable, nondecreasing, and bounded between zero and one.
We extend it by the constants zero and one outside this interval.

The corresponding DOS approximation is the almost-everywhere derivative of the
extended CDOS interpolant. On each interval $[\bar\theta_j,\bar\theta_{j+1})$,
$j=0,1,\ldots,M$, we use the polynomial representative
\begin{equation}\label{eq:DOS_CS}
\begin{split}
    \wt d_{cs}(t)
    :={}&3\left(
        \frac{d_j+d_{j+1}-2\Delta_j}{h_j^2}
        \right)(t-\bar\theta_j)^2 \\
      &+2\left(
        \frac{-2d_j-d_{j+1}+3\Delta_j}{h_j}
        \right)(t-\bar\theta_j)+d_j.
\end{split}
\end{equation}
Thus $\wt d_{cs}=\wt c_{cs}'$ almost everywhere. Consistent with the constant extension of $\wt c_{cs}$, we set
\begin{equation}\label{eq:DOS_CS_exterior}
    \wt d_{cs}(t)=0,
    \qquad
    t<\bar\theta_0
    \quad\text{or}\quad
    t\geq\bar\theta_{M+1}.
\end{equation}
% Here the $d_j$ are the nodal slopes used to construct the CDOS spline,
% whereas $\wt d_{cs}(t)$ denotes the resulting DOS function.  
Because
$\wt c_{cs}$ is nondecreasing and changes from zero to one, its derivative
is nonnegative almost everywhere and satisfies
\begin{equation}
    \int_{\mathbb R}\wt d_{cs}(t)\,dt=1.
\end{equation}
Thus the spline-derived DOS is a normalized, piecewise-quadratic density and does not require a Gaussian bandwidth in its construction. It is continuous at the interior knots, but may jump to zero at the two artificial endpoints. Values assigned to the DOS at those endpoints do not affect its integrals; the CDOS itself remains continuous, so its distributional derivative has no endpoint Dirac masses. The knot spacing and endpoint choices nevertheless impose an intrinsic resolution, and nonnegativity alone does not guarantee accuracy.

\subsection{Gaussian-Process Regression}\label{sec:GP}

An alternative to exact interpolation is Gaussian process (GP) regression, which places a probability distribution over functions and updates that distribution using observed data \cite{rasmussen2006gpml}. In our case, we aim to construct a smooth approximation to the CDOS from the midpoint data $(\bar\theta_j,\bar\vartheta_j)$, $j=1,\ldots,M$, augmented with the artificial endpoints $(\bar\theta_0,0)$ and $(\bar\theta_{M+1},1)$. These points will be referred to as the training points. The GP approach provides both a fitted function and a probabilistic description of uncertainty about that function.

A Gaussian process is a collection of random variables such that the function values at any finite set of inputs have a joint Gaussian distribution. We model the latent CDOS function $f$ by
\begin{equation}
    f\sim\mathcal{GP}(m,k),
    \qquad m(\theta)=\mathbb E[f(\theta)],
    \qquad k(\theta,\theta')=\operatorname{Cov}(f(\theta),f(\theta')).
\end{equation}
The mean function describes the prior trend, while the covariance kernel determines how function values at different inputs are related and controls the regularity of the process. To distinguish the latent function from the supplied training values, we use the Gaussian observation model
\begin{equation}
    \bar\vartheta_j=f(\bar\theta_j)+\varepsilon_j,
    \qquad \varepsilon_j\overset{\mathrm{iid}}{\sim}\mathcal N(0,\sigma_n^2),
    \qquad j=0,\ldots,M+1,
\end{equation}
where $\bar\vartheta_0=0$, $\bar\vartheta_{M+1}=1$, and the errors are independent of $f$. Here $\sigma_n^2$ models discrepancies between the training values and the latent function; it is not a Gaussian DOS-broadening parameter. With this Gaussian likelihood and fixed model hyperparameters, conditioning on the training values gives the posterior process.

%we place a Gaussian-process (GP) prior on a latent smooth CDOS approximation $f$:
%\begin{equation}\label{eq:gp_prior}
%    f\sim\mathcal{GP}\bigl(m(\cdot),k(\cdot,\cdot)\bigr),
%\end{equation}
%where $m$ is the prior mean and $k$ is a positive-definite covariance kernel.
%The training inputs and observations are
%\begin{equation}\label{eq:gp_data}
%    x=(\bar\theta_0,\bar\theta_1,\ldots,\bar\theta_M,\bar\theta_{M+1})^T,
%   \qquad
%    y=(0,\bar\vartheta_1,\ldots,\bar\vartheta_M,1)^T.
%\end{equation}
%We allow an observation-noise variance $\sigma_n^2$, so that
%$y_i=f(x_i)+\varepsilon_i$ with independent
%$\varepsilon_i\sim\mathcal N(0,\sigma_n^2)$.

Let $K$ be the covariance matrix with entries $K_{ij}=k(\bar\theta_i,\bar\theta_j)$ for $i,j \in \{0, \dots, M+1\}$, let
$k_*(\theta)=(k(\theta,\bar\theta_0),\ldots,k(\theta,\bar\theta_{M+1}))^T$, and let
$\bar{\boldsymbol{\theta}}=(\bar\theta_0,\ldots,\bar\theta_{M+1})^T$ and
$\bar{\boldsymbol{\vartheta}}=(0,\bar\vartheta_1,\ldots,\bar\vartheta_M,1)^T$ denote the training input and value vectors.
The vector $m(\bar{\boldsymbol{\theta}})$ contains the prior means evaluated at these inputs.
The posterior mean, which we use as the GP CDOS estimate, is
\begin{equation}\label{eq:gpmean}
\wt c_{gp}(\theta)
=m(\theta)+k_*(\theta)^T(K+\sigma_n^2I)^{-1}(\bar{\boldsymbol{\vartheta}}-m(\bar{\boldsymbol{\theta}})).
\end{equation}
For test points $\theta$ and $\theta'$, the posterior covariance of the latent function is
\begin{equation}\label{eq:gpcov}
k_{\mathrm{post}}(\theta,\theta')
=k(\theta,\theta')-k_*(\theta)^T(K+\sigma_n^2I)^{-1}k_*(\theta').
\end{equation}
In particular, a pointwise 95\% credible interval is
$\wt c_{gp}(\theta)\pm1.96\sqrt{k_{\mathrm{post}}(\theta,\theta)}$.
This interval describes uncertainty about the latent CDOS at a fixed input, conditional on the training data and hyperparameters. Because $\bar\vartheta_j$ is generally not the exact value of CDOS at $\bar\theta_j$, the posterior uncertainty is only meaningful in a relative sense.  Furthermore, the independent Gaussian observation errors are a modeling assumption, not a consequence of stochastic Lanczos sampling; midpoint errors share quadrature weights and nodes and can be correlated and biased. The posterior uncertainty therefore does not, by itself, quantify the full error relative to the exact CDOS or DOS. Nonetheless, the uncertainty quantification can be used to gauge, for example, which regions of the spectrum are more or less well resolved by the Lanczos quadrature. It can also be used to guide adaptive sampling, for example, by adding more Lanczos runs if some regions exhibit high posterior uncertainty.

When $\sigma_n^2=0$ and $K$ is nonsingular, the GP posterior mean interpolates the training data in exact arithmetic.
For $\sigma_n^2>0$, it smooths through the observations and need not pass exactly through them. Since $\bar\vartheta_j$ is generally not the exact value of CDOS at $\bar\theta_j$, it is preferable to allow $\sigma_n^2>0$ and optimize it along with the kernel hyperparameters. 

The prior mean $m(\theta)$ is often set to zero when no specific trend is prescribed. A common covariance kernel is the squared exponential,
\begin{equation}
k(\theta, \theta') = \sigma_f^2 \exp\left( - \frac{(\theta - \theta')^2}{2\ell^2} \right),
\end{equation}
where $\sigma_f^2$ controls the prior variance and $\ell$ is the correlation length. Mat\'ern kernels provide another family with an adjustable degree of regularity. The kernel hyperparameters and observation variance can be estimated from the training data, for example by maximizing the marginal likelihood. Kernel choice determines the smoothness assumptions; Gaussianity alone does not ensure smooth sample functions.

To accommodate variation on two correlation scales, we also consider a double-Gaussian covariance,
\begin{equation}\label{eq:gp_dual_kernel}
    k(\theta,\theta')=k_1(\theta,\theta')+k_2(\theta,\theta'),
    \qquad
    k_r(\theta,\theta')=\sigma_{f,r}^2
    \exp\!\left(-\frac{(\theta-\theta')^2}{2\ell_r^2}\right),
    \quad r=1,2.
\end{equation}
Each component has its own variance and correlation length. If the fitted lengths differ, the longer-scale component can represent a broad trend while the shorter-scale component captures finer variation. This kernel corresponds to the sum of two independent Gaussian processes, not a mixture of two Gaussian distributions. For this model we take $m(\theta)=0$; the posterior formulas~\eqref{eq:gpmean}--\eqref{eq:gpcov} remain unchanged, with $K$ and $k_*$ formed from the summed kernel.

Writing $C=K+\sigma_n^2I$ and
$\mathbf r=\bar{\boldsymbol{\vartheta}}-m(\bar{\boldsymbol{\theta}})$,
the negative log marginal likelihood is
\begin{equation}\label{eq:gp_nlml}
    \mathcal L(\boldsymbol\eta)
    =\frac12\mathbf r^TC^{-1}\mathbf r
     +\frac12\log\det C+\frac{M+2}{2}\log(2\pi).
\end{equation}
For the double-Gaussian model,
$\boldsymbol\eta=(\sigma_{f,1}^2,\ell_1,\sigma_{f,2}^2,\ell_2,\sigma_n^2)$
contains the five positive parameters to be fitted by minimizing $\mathcal L$.

\begin{figure}[htbp]
\centering
\includegraphics[width=\textwidth]{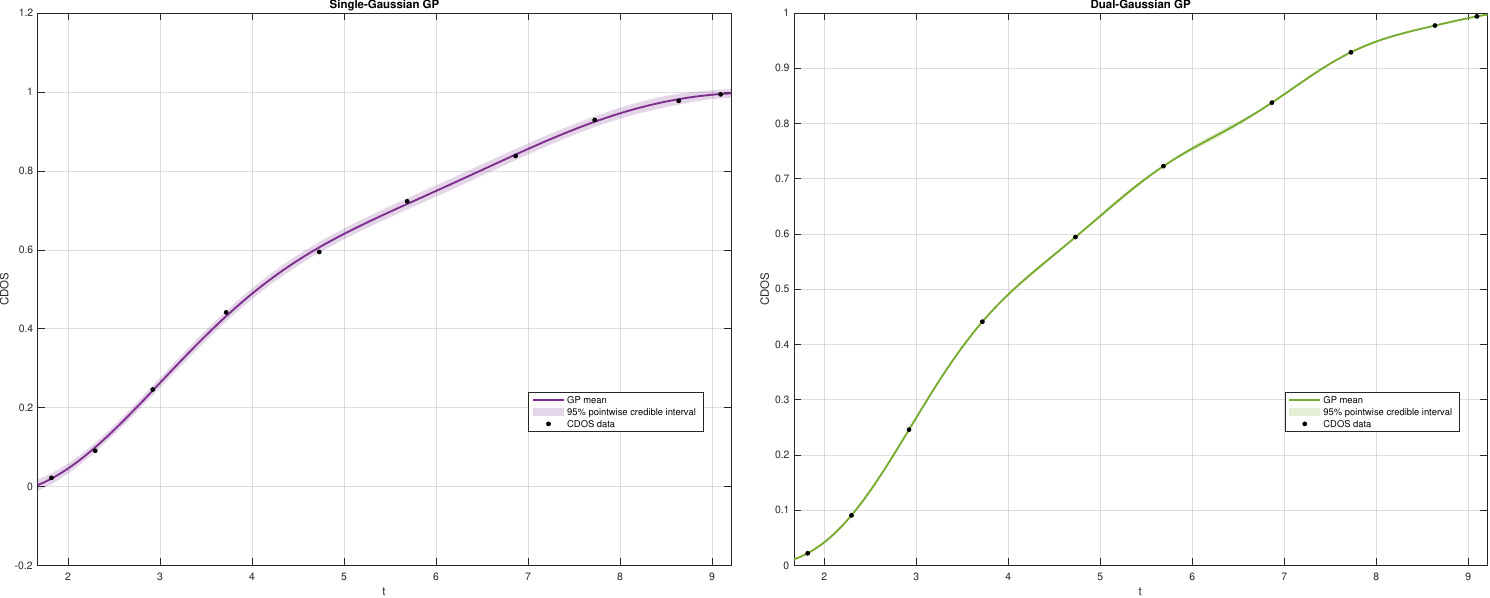}
\caption{GP CDOS approximations for \texttt{mesh1e1}, fitted to the averaged midpoint data from the right panel of Figure~\ref{fig:mesh1e1-CDOS-Lanczos}. Left: a single squared-exponential kernel. Right: the sum of two squared-exponential kernels~\eqref{eq:gp_dual_kernel}. Black dots show the midpoint data, colored curves show posterior means, and shaded regions show pointwise 95\% credible intervals for the latent CDOS, conditional on the fitted hyperparameters. The double-Gaussian band is barely visible at this scale.}
\label{fig:mesh1e1-CDOS-single-double-GP}
\end{figure}

Figure~\ref{fig:mesh1e1-CDOS-single-double-GP} compares the single-Gaussian and double-Gaussian GP fits to the same averaged midpoint data shown in the right panel of Figure~\ref{fig:mesh1e1-CDOS-Lanczos}. In this example, the single-Gaussian posterior mean does not pass exactly through all training points, and its pointwise 95\% credible band has a modest but visible width. The double-Gaussian posterior mean follows the training points more closely, while its band is barely visible at the plotted scale. These are properties of the fitted models in this example, not general interpolation or uncertainty guarantees associated with the number of kernel components.

Closer agreement with the training points does not, by itself, establish greater accuracy relative to the exact CDOS. The midpoint data are themselves approximations, affected by stochastic probing and finite-step quadrature errors. Likewise, a narrower band reflects the posterior uncertainty under the fitted covariance and observation-noise model, conditional on the fitted hyperparameters; it does not imply that these errors or hyperparameter uncertainty have been accounted for. In particular, neither closer agreement with the CDOS training points nor a narrower CDOS band guarantees a more accurate DOS derivative.

The smooth covariance kernels above also allow the DOS to be obtained by analytic differentiation. With a differentiable mean function, the derivative $f'$ is a Gaussian process whose posterior mean and pointwise variance are
\begin{equation}\label{eq:gp_dos}
\begin{aligned}
    \wt d_{gp}(\theta)
    &:=\frac{d}{d\theta}\wt c_{gp}(\theta)
      =m'(\theta)+h_*(\theta)^TC^{-1}\mathbf r,\\
    v_{d,gp}(\theta)
    &=\left.\frac{\partial^2 k(\theta,\theta')}
                   {\partial\theta\,\partial\theta'}\right|_{\theta'=\theta}
      -h_*(\theta)^TC^{-1}h_*(\theta),
\end{aligned}
\end{equation}
where $h_*(\theta)=\frac{d}{d\theta}k_*(\theta)$. For the double-Gaussian kernel, these derivatives are explicitly
\begin{equation}\label{eq:gp_dual_derivatives}
\begin{aligned}
    {[h_*(\theta)]_j}
    &=-\sum_{r=1}^2\frac{\theta-\bar\theta_j}{\ell_r^2}
                       k_r(\theta,\bar\theta_j),
      \qquad j=0,\ldots,M+1,\\
    \left.\frac{\partial^2 k(\theta,\theta')}
                   {\partial\theta\,\partial\theta'}\right|_{\theta'=\theta}
    &=\frac{\sigma_{f,1}^2}{\ell_1^2}
      +\frac{\sigma_{f,2}^2}{\ell_2^2}.
\end{aligned}
\end{equation}
The single-Gaussian expressions follow by retaining only one kernel component. A pointwise 95\% credible interval for the latent derivative is
$\wt d_{gp}(\theta)\pm1.96\sqrt{v_{d,gp}(\theta)}$,
conditional on the fitted hyperparameters. This variance concerns the derivative of the latent CDOS, not noisy observations, so no observation-noise variance is added to it.

The standard unconstrained GP model described here does not automatically enforce monotonicity or the range $[0,1]$. Consequently, its DOS derivative need not be nonnegative. Moreover, its integral over the training interval equals $\wt c_{gp}(\bar\theta_{M+1})-\wt c_{gp}(\bar\theta_0)$, which need not be exactly one when the endpoint observations are fitted with nonzero noise variance. Such shape and endpoint constraints require additional modeling assumptions, whereas they are built into the monotone spline construction above.

Monotonicity information can, for example, be introduced through likelihood factors favoring positive derivatives at auxiliary locations \cite{riihimaki2010monotonicity}. Such finite, soft constraints do not guarantee global monotonicity. They also lead to a non-Gaussian posterior, so the exact Gaussian-conditioning formulas above must be replaced by an appropriate inference approximation. On the full real line, a zero-mean squared-exponential posterior mean tends to zero at both ends and is not a normalized CDOS; any finite-domain extension used for DOS validation must therefore be specified separately.

\section{Metrics for Evaluating DOS Approximations}\label{sec:metrics}

The exact DOS~\eqref{eq:dos} is a sum of Dirac distributions, so it cannot be compared pointwise with a function-valued approximation using ordinary function norms. We therefore evaluate DOS approximations at a specified finite spectral resolution. Four measures are considered: relative $L^\infty$ error, relative $L^2$ error, Jensen--Shannon divergence, and cosine error. The first two assess discrepancies in DOS amplitudes, whereas the latter two compare normalized curves and are therefore insensitive to a uniform positive rescaling of either curve.

\subsection{Reference Curves and Validation Resolution}

Fix a validation bandwidth $\sigma>0$ and an evaluation interval $I=[a,b]$. The reference curve is the exact Gaussian-regularized DOS
\begin{equation}\label{eq:metrics_reference}
    r(t):=d_\sigma(t)=\frac1n\sum_{j=1}^n g_\sigma(t-\lambda_j).
\end{equation}
For a common-resolution comparison, the approximate spectral information is convolved with the same Gaussian. In particular, a function-valued estimate $\widehat d$, such as the derivative of a spline or GP CDOS approximation, gives
\begin{equation}\label{eq:metrics_compared_dos}
    s(t):=(g_\sigma*\widehat d)(t)
    =\int_{\mathbb R}g_\sigma(t-u)\widehat d(u)\,du.
\end{equation}
For discrete Lanczos quadrature data, the corresponding curve is the weighted Gaussian sum $s(t)=\widehat d_{\sigma,M}(t)$. This comparison measures agreement of the underlying spectral information after applying a common resolution filter. The validation bandwidth is an evaluation parameter, not a parameter required to construct the spline or GP approximation.

Alternatively, one can compare $s(t)=\widehat d(t)$ directly with $r(t)=d_\sigma(t)$. This asks how closely the returned DOS curve approximates a prescribed Gaussian target, rather than comparing both underlying representations after common smoothing. The formulas below apply to either convention, but the convention and $\sigma$ must be stated when reporting results. In particular, applying another Gaussian convolution to an already Gaussian-broadened curve changes its effective resolution.

All norms and integrals below are taken over $I$, which should include the spectral region of interest and sufficient Gaussian tails. The compared curves are not independently normalized before computing the relative norm errors or cosine error. For Jensen--Shannon divergence, a separate probability normalization is required. Thus the choice of interval matters, especially if appreciable spectral mass lies outside it.

\subsection{Relative Maximum Error}

The relative $L^\infty$ error is
\begin{equation}\label{eq:metrics_relative_linf}
    E_{\infty,\mathrm{rel}}
    :=\frac{\|s-r\|_{L^\infty(I)}}{\|r\|_{L^\infty(I)}}
    =\frac{\sup_{t\in I}|s(t)-r(t)|}{\sup_{t\in I}|r(t)|}.
\end{equation}
For the continuous curves considered here, this is the largest pointwise discrepancy divided by the maximum reference height. It detects a badly misplaced or incorrectly sized peak, a spurious narrow spike, or a deep local trough even if the error is confined to a small interval. For example, $E_{\infty,\mathrm{rel}}=0.1$ means that the largest absolute discrepancy is one tenth of the maximum reference height. This is a global normalization, not the maximum of $|s(t)-r(t)|/r(t)$, which can be dominated by tiny reference values in spectral gaps or tails.

The metric is useful when the worst local error is important, but it does not indicate how widely that error is distributed across the spectrum. A narrow spike and a broad discrepancy with the same maximum amplitude can have identical $L^\infty$ errors. The relative form is dimensionless and can exceed one; it remains well-defined for a signed approximation $s$ whenever the reference norm is nonzero.

\subsection{Relative \texorpdfstring{$L^2$}{L2} Error}

The relative $L^2$ error is
\begin{equation}\label{eq:metrics_relative_l2}
    E_{2,\mathrm{rel}}
    :=\frac{\|s-r\|_{L^2(I)}}{\|r\|_{L^2(I)}}
    =\left(\frac{\int_I|s(t)-r(t)|^2\,dt}
                       {\int_I|r(t)|^2\,dt}\right)^{1/2}.
\end{equation}
Unlike the maximum error, this measure accounts for both the amplitude and the spectral extent of the discrepancy. An error of approximately constant height $h$ over an interval of width $w$ contributes approximately $h^2w$ to the squared numerator. Consequently, a broad mismatch or persistent oscillations can produce a substantial $L^2$ error even when no single discrepancy is extreme. Squaring the error gives greater weight to large absolute deviations, so inaccuracies near prominent peaks can dominate smaller discrepancies elsewhere.

This metric provides an overall measure of curve accuracy and complements the local worst-case information supplied by $E_{\infty,\mathrm{rel}}$. It is dimensionless, is not bounded by one, and is valid for signed approximations. Both relative norm errors retain sensitivity to an incorrect overall DOS amplitude or spectral mass; neither requires the approximate curve to be a probability density.

\subsection{Jensen--Shannon Divergence}

To compare how spectral mass is distributed, suppose that $r$ and $s$ are nonnegative on $I$ and have positive integrals. Define probability densities on this interval by
\begin{equation}\label{eq:metrics_probability_normalization}
    p(t):=\frac{r(t)}{\int_I r(u)\,du},\qquad
    q(t):=\frac{s(t)}{\int_I s(u)\,du},\qquad
    m_{pq}(t):=\frac{p(t)+q(t)}2.
\end{equation}
The Jensen--Shannon divergence is the average relative entropy of each density with respect to their equal mixture~\cite{lin1991divergence}:
\begin{equation}\label{eq:metrics_js}
\begin{aligned}
    D_{\mathrm{JS}}(p,q)
    := {}&\frac12\int_I p(t)\log\!\left(\frac{p(t)}{m_{pq}(t)}\right)\,dt\\
         &+\frac12\int_I q(t)\log\!\left(\frac{q(t)}{m_{pq}(t)}\right)\,dt.
\end{aligned}
\end{equation}
We use natural logarithms and define each integrand to be zero wherever its prefactor is zero, including where both densities vanish. The divergence is symmetric and satisfies $0\leq D_{\mathrm{JS}}\leq\log 2$, with zero precisely when $p=q$ almost everywhere. It remains finite when one density vanishes where the other is positive. Using base-two logarithms instead gives a value in $[0,1]$. The square root is the Jensen--Shannon distance~\cite{endres2003metric}; here we report the divergence itself.

This measure assesses the distinguishability of the two normalized spectral distributions. It uses density ratios weighted by probability mass, rather than squared absolute height differences as in $L^2$. Separate normalization makes it insensitive to an overall positive multiplicative error: if $s=\gamma r$ for $\gamma>0$, then $D_{\mathrm{JS}}=0$, even when the unnormalized masses disagree. It therefore complements, but does not replace, the relative norm errors.

Nonnegativity is essential. If a GP CDOS approximation is not monotone, its derivative can be negative and hence fail to define a probability density. Whenever the compared curve $s$ has a genuinely negative region, the Jensen--Shannon divergence above is undefined. Replacing $s$ by its positive part $\max(s,0)$ and renormalizing would evaluate a modified approximation and can conceal this failure. Such a projected value must be identified separately, rather than substituted for the divergence of the original curve. Roundoff-level negative values may be treated as zero under a stated numerical tolerance. Gaussian validation smoothing can remove some negative regions, so a valid divergence for a smoothed curve does not establish that the original DOS was nonnegative.

\subsection{Cosine Error}

For nonzero $L^2$ curves, define the cosine similarity and cosine error by
\begin{equation}\label{eq:metrics_cosine}
    S_{\mathrm{cos}}(r,s)
    :=\frac{\int_I r(t)s(t)\,dt}
            {\|r\|_{L^2(I)}\|s\|_{L^2(I)}},
    \qquad E_{\mathrm{cos}}:=1-S_{\mathrm{cos}}(r,s).
\end{equation}
The similarity is the normalized inner product of the curves, viewed as vectors in $L^2(I)$. It measures their alignment without subtracting their means. Equivalently,
\begin{equation}\label{eq:metrics_cosine_normalized}
    E_{\mathrm{cos}}
    =\frac12\left\|
       \frac{r}{\|r\|_{L^2(I)}}-\frac{s}{\|s\|_{L^2(I)}}
      \right\|_{L^2(I)}^2.
\end{equation}
Thus cosine error compares shape after $L^2$ normalization, whereas Jensen--Shannon divergence compares probability densities after mass normalization. Like squared-error measures, cosine error can be strongly influenced by high-amplitude peaks. It is zero when the curves differ only by a positive scale factor, so it does not detect an overall mass error.

For nonnegative curves, $0\leq E_{\mathrm{cos}}\leq1$, with one corresponding to zero overlap. For a signed approximation it remains algebraically defined and lies in $[0,2]$, provided both norms are nonzero. This makes it available even when Jensen--Shannon divergence is undefined, but a small cosine error is not a certificate of DOS nonnegativity. If either curve has zero norm, cosine error is undefined. Although often called cosine distance, $1-S_{\mathrm{cos}}$ is not a metric in the strict triangle-inequality sense.

The distinction from relative $L^2$ error can be made explicit. With $\eta=\|s\|_{L^2(I)}/\|r\|_{L^2(I)}$, direct expansion gives
\begin{equation}\label{eq:metrics_l2_cosine_relation}
    E_{2,\mathrm{rel}}^2=(\eta-1)^2+2\eta E_{\mathrm{cos}}.
\end{equation}
The relative $L^2$ error therefore combines an overall norm mismatch with a shape mismatch. For $s=\gamma r$ with $\gamma>0$, both relative norm errors equal $|\gamma-1|$, while the cosine error and Jensen--Shannon divergence vanish.

\subsection{Evaluation on a Discrete Grid}

Let $a=t_1<\cdots<t_N=b$ be a common evaluation grid and let $r_i=r(t_i)$ and $s_i=s(t_i)$. To approximate the integrals consistently on a possibly nonuniform grid, we use trapezoidal weights
\begin{equation}\label{eq:metrics_quadrature_weights}
    w_1=\frac{t_2-t_1}{2},\qquad
    w_i=\frac{t_{i+1}-t_{i-1}}2\quad(2\leq i\leq N-1),\qquad
    w_N=\frac{t_N-t_{N-1}}2.
\end{equation}
The relative errors and cosine similarity are then approximated by
\begin{equation}\label{eq:metrics_discrete}
\begin{aligned}
    E_{\infty,\mathrm{rel}}&\approx
       \frac{\max_i|s_i-r_i|}{\max_i|r_i|},\qquad
    E_{2,\mathrm{rel}}\approx
       \left(\frac{\sum_i w_i(s_i-r_i)^2}{\sum_i w_i r_i^2}\right)^{1/2},\\
    S_{\mathrm{cos}}&\approx
       \frac{\sum_i w_i r_i s_i}
            {\sqrt{\sum_i w_i r_i^2}\sqrt{\sum_i w_i s_i^2}}.
\end{aligned}
\end{equation}
For Jensen--Shannon divergence, normalize using these same weights $p_i=r_i/\sum_jw_jr_j$ and $q_i=s_i/\sum_jw_js_j$, and replace the integrals in~\eqref{eq:metrics_js} by weighted sums. Merely normalizing by the sum of sampled heights would not represent the intended probability densities on a nonuniform grid.

Grid refinement is needed to resolve narrow peaks and oscillations; in particular, the sampled maximum can miss a discrepancy between grid points. The interval and the treatment of the approximation outside the sampled domain must also be specified when performing validation convolution. If the sampled DOS is extended by zero, the domain must be wide enough to avoid truncating relevant mass. All four measures should be compared at the same declared resolution and over the same interval. Taken together, they distinguish worst-case local error, integrated amplitude error, probability-distribution mismatch, and scale-independent $L^2$ shape mismatch.

\section{Numerical Results}\label{sec:results}

We consider six real symmetric matrices from electronic-structure calculations, interacting quantum many-body models, nuclear configuration interaction, and spatial discretizations. This selection provides variation not only in application and matrix dimension, but also in spectral range, clustering, gaps, and the location of high-density regions. The full eigenvalue sets are used to construct reference spectra and Gaussian-regularized DOS curves for validation; the practical DOS estimators do not require full diagonalization. The ideal-spline reference deliberately uses spectral information.

\subsection{Test Matrices}\label{sec:test_matrices}

Table~\ref{tab:six_test_matrices} summarizes the six problems and the Gaussian widths used to visualize their exact DOS. The matrices \texttt{benzene}, \texttt{nd3k}, and \texttt{fv1} are drawn from the SuiteSparse Matrix Collection~\cite{davis2011university}; the remaining three are supplied many-body Hamiltonians. All dimensions refer to the actual matrices tested, including the restricted many-body sectors described below.

\begin{table}[tbp]
\centering
\small
\caption{Test matrices, application areas, and Gaussian widths for the exact DOS curves in Figure~\ref{fig:six_matrix_exact_DOS}. The widths are expressed in each matrix's spectral units.}
\label{tab:six_test_matrices}
\begin{tabular}{l r p{0.43\textwidth} c}
\hline
Matrix & $n$ & Application & $\sigma$ \\
\hline
\texttt{benzene} & $8{,}219$ & Electronic structure: Kohn--Sham Hamiltonian & $0.3$ \\
\texttt{HubbardL8} & $3{,}136$ & Interacting fermions on a periodic chain & $0.2$ \\
\texttt{HeisenbergL16} & $12{,}870$ & Disordered quantum spin chain & $0.1$ \\
\texttt{Li7Nmax2} & $1{,}961$ & Nuclear configuration-interaction Hamiltonian & $0.5$ \\
\texttt{nd3k} & $9{,}000$ & Molecular normal mode analysis problem & $0.3$ \\
\texttt{fv1} & $9{,}604$ & Finite-element discretization of the Laplace equation & $0.1$ \\
\hline
\end{tabular}
\end{table}

\paragraph{Benzene: electronic structure.}
The \texttt{benzene} matrix is a real symmetric Kohn--Sham Hamiltonian from a density-functional-theory calculation using the PARSEC real-space pseudopotential method~\cite{cts:94,che-PDFM00}. It has dimension $n=8219$ and $242{,}669$ nonzero entries. Its eigenvalues represent electronic energy levels in the discretized problem. This example connects DOS estimation directly to its electronic-structure applications and provides an asymmetric spectral density with a broad dominant peak and a higher-energy tail.

\paragraph{HubbardL8: interacting fermions.}
This matrix represents the one-dimensional spin-$\tfrac12$ Hubbard model on $L=8$ sites with periodic boundary conditions. The Hamiltonian is
\begin{equation}\label{eq:test_hubbard}
    H_{\mathrm{Hub}}
    =-t_{\mathrm{hop}}\sum_{\langle i,j\rangle}\sum_{s\in\{\uparrow,\downarrow\}}
       \left(c_{i,s}^{\dagger}c_{j,s}+c_{j,s}^{\dagger}c_{i,s}\right)
      +U\sum_{i=1}^{L}n_{i,\uparrow}n_{i,\downarrow},
\end{equation}
where $c_{i,s}^{\dagger}$ and $c_{i,s}$ create and annihilate a fermion of spin $s$ at site $i$, $n_{i,s}=c_{i,s}^{\dagger}c_{i,s}$, and $\langle i,j\rangle$ denotes nearest-neighbor bonds, including the bond between sites $L$ and $1$. We use hopping strength $t_{\mathrm{hop}}=1$ and on-site repulsion $U=4$, with no additional on-site potential or chemical-potential term. Restricting to $N_\uparrow=N_\downarrow=3$ particles gives a sparse real symmetric matrix of dimension
\[
    n=\binom{8}{3}\binom{8}{3}=3136.
\]
Unlike the one-particle electronic-structure example, its eigenvalues describe energies of interacting many-particle configurations. Its regularized DOS exhibits several resolved local peaks within a broad central concentration.

\paragraph{HeisenbergL16: a disordered spin chain.}
We use the nearest-neighbor spin-$\tfrac12$ Heisenberg Hamiltonian on an open chain of $L=16$ sites,
\begin{equation}\label{eq:test_heisenberg}
    H_{\mathrm{Heis}}
    =\sum_{i=1}^{L-1}\left(S_i^xS_{i+1}^x+S_i^yS_{i+1}^y+S_i^zS_{i+1}^z\right)
      +\sum_{i=1}^{L}h_iS_i^z,
\end{equation}
where $S_i^\alpha$ denotes the spin operator at site $i$ in direction $\alpha$, and the exchange coupling is set to one. The on-site fields are sampled independently and uniformly from $[-w,w]$ with disorder strength $w=5$. The Hamiltonian preserves total magnetization. We take its largest fixed-magnetization block, containing eight spins up and eight spins down, so that
\[
    n=\binom{16}{8}=12870.
\]
This sparse real symmetric problem complements the periodic Hubbard model through its open boundaries, spin interactions, and spatial disorder. Its DOS has a broad, approximately bell-shaped envelope with finer local structure.

\paragraph{Li7Nmax2: nuclear configuration interaction.}
The \texttt{Li7Nmax2} matrix is a configuration-interaction Hamiltonian for the nucleus ${}^{7}\mathrm{Li}$ with truncation parameter $N_{\max}=2$, which results in dimension $n=1961$. The finite configuration space represents many-nucleon states, and the eigenvalues are the corresponding energy levels within that truncation. This provides a nuclear many-body example distinct from both electronic-structure and lattice-model Hamiltonians. Among the displayed curves, it has a broad spectral range and substantial resolved local variation in the DOS.

\paragraph{nd3k: molecular normal mode analysis problem.}
The matrix \texttt{nd3k} belongs to the ND problem set in SuiteSparse. The matrix originates from ~\cite{tavgnca}. It is the Hessian of the potential energy surface associated with a 3000-atom polymer. It is real symmetric positive definite, has dimension $n=9000$, and contains $3{,}279{,}690$ nonzero entries. Relative to the other SuiteSparse examples in this test set, its matrix is substantially more densely populated. Its spectrum contains a pronounced low-eigenvalue cluster separated by a large gap from a higher-eigenvalue group, making it useful for examining how DOS approximations represent strongly nonuniform spectral distributions.

\paragraph{fv1: an elliptic partial differential equation.}
The \texttt{fv1} matrix comes from a finite-element discretization of the Laplace equation on a two-dimensional $100\times100$ mesh in the Norris collection. The resulting matrix has dimension $n=9604$ and $85{,}264$ nonzero entries and is real symmetric positive definite. Its eigenvalues occupy a much narrower numerical range than those of \texttt{nd3k}, and its DOS is strongly concentrated near the upper end of the spectrum. Thus the two positive-definite discretization problems present quite different density profiles despite sharing broad structural properties.

\subsection{Spectral Distributions and Exact DOS}\label{sec:test_spectra}

\begin{figure}[tbp]
\centering
\includegraphics[width=\textwidth]{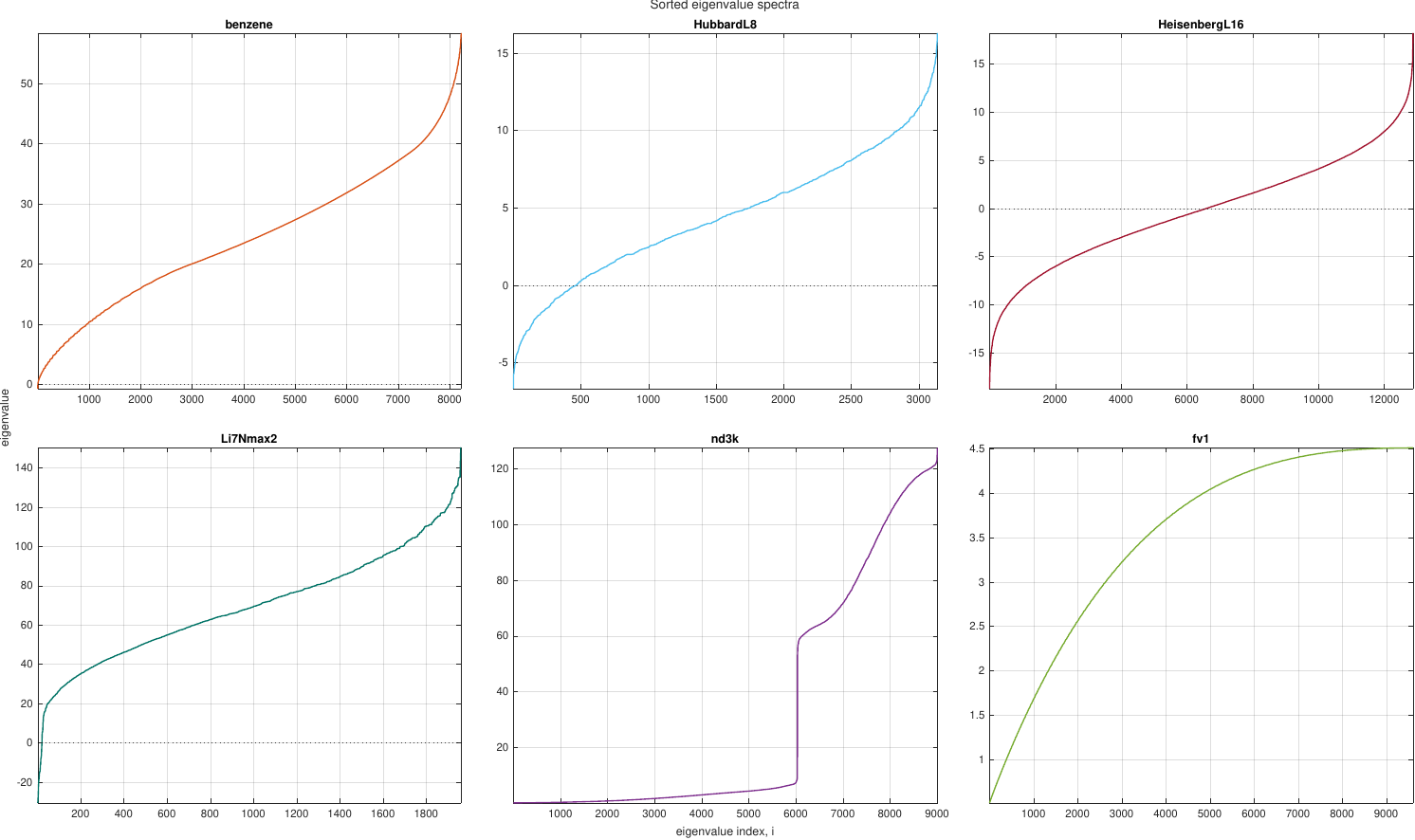}
\caption{Sorted eigenvalues $\lambda_i$ versus index $i$ for the six test matrices. The top row shows \texttt{benzene}, \texttt{HubbardL8}, and \texttt{HeisenbergL16}; the bottom row shows \texttt{Li7Nmax2}, \texttt{nd3k}, and \texttt{fv1}. Flat portions indicate eigenvalue concentrations, whereas steep portions indicate sparse spectral regions or gaps.}
\label{fig:six_matrix_eigenvalue_spectra}
\end{figure}

Figure~\ref{fig:six_matrix_eigenvalue_spectra} plots each sorted eigenvalue $\lambda_i$ against its index $i$. These plots show the spectral ranges and large gaps, including the separation between the two main eigenvalue groups in \texttt{nd3k}. However, local concentrations are encoded indirectly through changes in the slope: slowly increasing segments contain many eigenvalues within a small spectral interval, while steep segments indicate sparse regions or gaps. Consequently, detailed density features are less immediately apparent than in a DOS plot.

\begin{figure}[tbp]
\centering
\includegraphics[width=\textwidth]{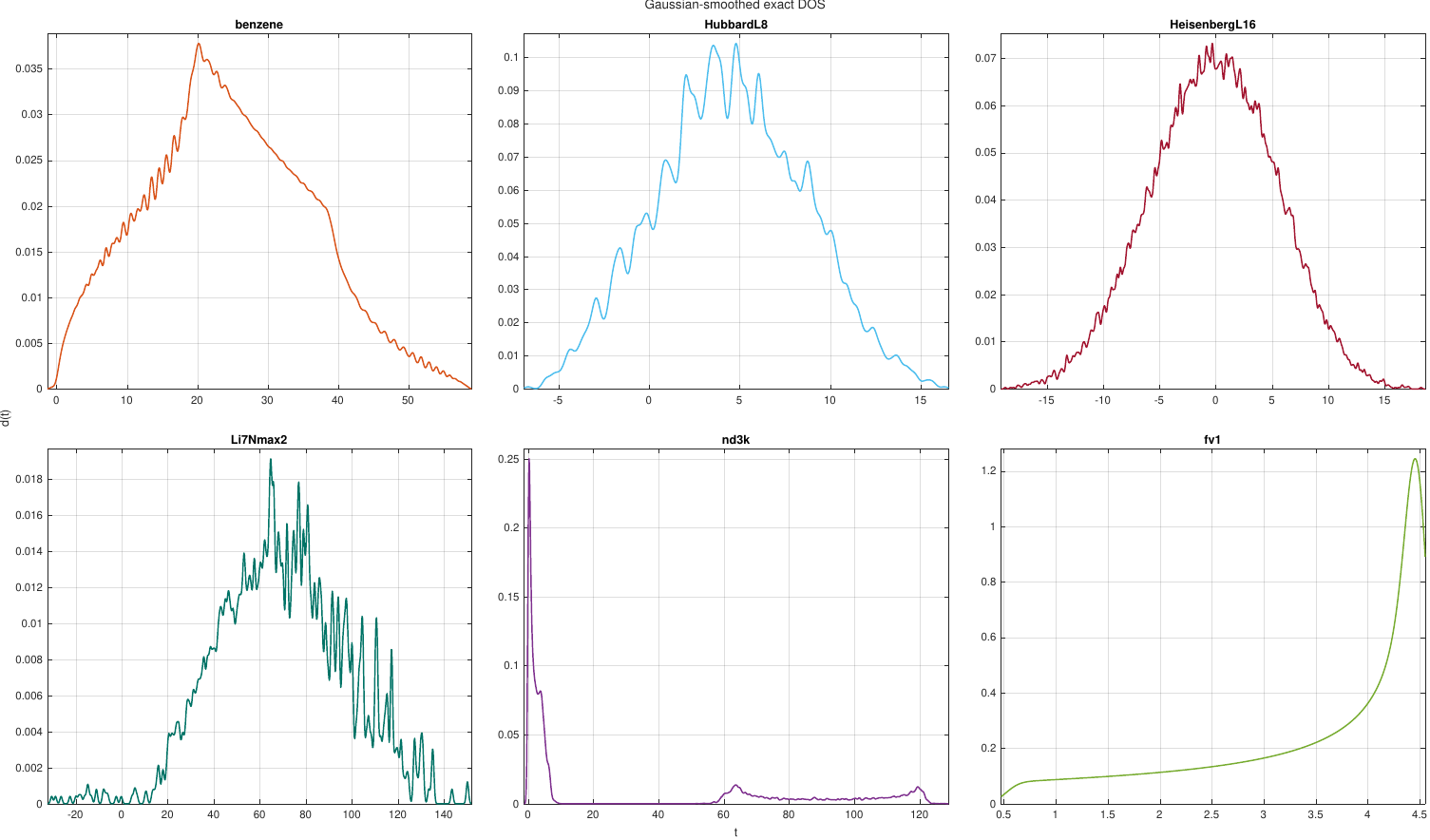}
\caption{Exact Gaussian-regularized DOS $d_\sigma(t)$ for the six test matrices, in the same panel order as Figure~\ref{fig:six_matrix_eigenvalue_spectra}. The top-row widths are $\sigma=0.3$, $0.2$, and $0.1$; the bottom-row widths are $\sigma=0.5$, $0.3$, and $0.1$. Each curve is the normalized Gaussian sum over the full spectrum. The individually chosen widths balance visible spectral detail with readability; axis scales differ between panels.}
\label{fig:six_matrix_exact_DOS}
\end{figure}

Figure~\ref{fig:six_matrix_exact_DOS} presents the same eigenvalue information through the exact Gaussian-regularized DOS $d_\sigma(t)$ in~\eqref{eq:regDOS}. Here ``exact'' means that the Gaussian sum uses all eigenvalues, rather than Ritz values or a fitted CDOS derivative. The widths are $0.3$, $0.2$, $0.1$, $0.5$, $0.3$, and $0.1$ for \texttt{benzene}, \texttt{HubbardL8}, \texttt{HeisenbergL16}, \texttt{Li7Nmax2}, \texttt{nd3k}, and \texttt{fv1}, respectively. They were selected by visual inspection to retain informative peaks and concentrations while limiting fine-scale oscillations that would obscure the overall distribution. They are not claimed to be optimal widths, and their numerical values must be interpreted relative to the different spectral scales.

The DOS plots make the differences between these problems immediately visible: the asymmetric profile of \texttt{benzene}, the multiple central peaks of \texttt{HubbardL8}, the broad central envelope of \texttt{HeisenbergL16}, the structured density of \texttt{Li7Nmax2}, the low-end concentration and separated spectral groups of \texttt{nd3k}, and the upper-end concentration of \texttt{fv1}. They illustrate the practical value of DOS visualization: regions containing many eigenvalues appear directly as peaks, rather than having to be inferred from the slope of a sorted-spectrum curve. Gaussian broadening does suppress detail below the selected resolution, so the two representations are complementary rather than interchangeable. Together, these six problems provide a varied test set spanning several application areas and qualitatively different spectral distributions.

\subsection{Polynomial Comparison Method}\label{sec:comparison_methods}
\subsubsection{Kernel Polynomial Method}

The kernel polynomial method (KPM) approximates the exact DOS of a matrix $A$ by expanding the Dirac delta distribution in a basis of Chebyshev polynomials $T_k(t)=\cos\bigl(k\arccos(t)\bigr)$.
For the formulas in this subsection, we use a scaled spectral coordinate so that the spectrum lies strictly inside $(-1,1)$. In general, one applies an affine transformation $A_s=(A-b_0I)/a_0$, with $a_0>0$, chosen from enclosing spectral bounds. The strict enclosure avoids placing eigenvalues at the singular endpoints of the Chebyshev weight.
For notational simplicity, $A$, $\lambda_j$, and $d$ in the following expansion formulas refer to the scaled matrix, its eigenvalues, and its DOS.

Following Silver and R\"oder~\cite{SilverRoder1994}, we expand the modified density $\hat{d}(t)=\sqrt{1-t^2}d(t)$ as
\begin{equation}\label{eq:exp1}
\hat d(t)  = \sum_{k=0}^\infty  \mu_k T_k(t).
\end{equation}
Eq.~\eqref{eq:exp1} should be understood in the sense of distributions, i.e. for any test function $g\in \mathcal{S}$,
\[
\int_{-1}^{1}\hat{d}(s) g(s) \ ds = \int_{-1}^{1}
\sum_{k=0}^{\infty} \mu_{k} T_{k}(s) g(s) \ ds.
\]
It can be shown that the coefficients of the expansion can be evaluated as
\begin{equation}
\mu_k = \frac{2-\delta_{k0}}{n\pi}  \sum_{j=1}^n  T_k(\lambda_j),
\label{eq:expmu}
\end{equation}
where $\delta_{k0}$ is the Kronecker delta (so $2-\delta_{k0}=1$ for $k=0$, and $2$ otherwise).

The coefficients can be estimated from stochastic traces. For example, with independent raw Gaussian probes $g^{(l)}\sim\mathcal N(0,I)$, define
\begin{equation}\label{eq:kpm_stochastic_coefficients}
    \widehat\mu_k
    =\frac{2-\delta_{k0}}{n\pi\nvec}
       \sum_{l=1}^{\nvec}(g^{(l)})^T T_k(A)g^{(l)}.
\end{equation}
Alternatively, normalized isotropic directions $v^{(l)}$ can be used with prefactor $(2-\delta_{k0})/(\pi\nvec)$; this convention fixes the estimated zeroth moment exactly. In either case, a common set of probes is used across all degrees. The recurrence
\[
    T_{k+1}(A)v=2A\,T_k(A)v-T_{k-1}(A)v
\]
requires one matrix--vector product per additional degree; each moment then requires an inner product. With estimated coefficients, the undamped DOS approximation is
\begin{equation}
\tilde d_M(t)
\;=\;
\frac{1}{\sqrt{1-t^2}}
\sum_{k=0}^M \widehat\mu_k \,T_k(t).
\label{eq:tdosFinal}
\end{equation}

Analogous expansions can be constructed with other orthogonal polynomial families, such as Legendre polynomials, using their corresponding weights and coefficients.
Moreover, it is possible to first regularize the DOS by replacing each $\delta(t-\lambda_j)$ in \eqref{eq:dos} with a Gaussian or Lorentzian function and then expand that smooth density.
The analytical expressions for the resulting expansion coefficients are derived in \cite{LSY}.

A notable caveat of the undamped truncation~\eqref{eq:tdosFinal} is that it need not remain nonnegative and can exhibit Gibbs oscillations. In the experiments below, we instead use the Jackson-damped variant, denoted \emph{KPM--Jackson},
\begin{equation}\label{eq:kpm_jackson_dos}
    \tilde d_M^{\mathrm J}(t)
    =\frac{1}{\sqrt{1-t^2}}\sum_{k=0}^{M}
       g_{k,M}^{\mathrm J}\widehat\mu_kT_k(t),
\end{equation}
where $g_{k,M}^{\mathrm J}$ are the Jackson damping coefficients, with $g_{0,M}^{\mathrm J}=1$. Jackson damping gives a positive reconstruction kernel and suppresses Gibbs oscillations, at the cost of broadening fine spectral features~\cite{kpmsurvey2006}. Consequently, nonnegativity is not an advantage of the spline over a correctly implemented KPM--Jackson estimator; the comparison concerns accuracy for a given number of matrix--vector products (matvecs). After computation in the scaled coordinate $x=(t-b_0)/a_0$, the physical-coordinate density is recovered as $a_0^{-1}\tilde d_M^{\mathrm J}((t-b_0)/a_0)$. The common Gaussian validation convolution described below is applied in addition to Jackson damping.

Chen's spectrum-adaptive KPM~\cite{chen2023adaptive} computes modified polynomial moments from the small Lanczos tridiagonal matrix via quadrature, allowing the spectral interval and polynomial reference density to be selected or revised after the matrix--vector products have been completed. This flexibility avoids repeating the large-matrix computation when these choices change, although the polynomial reconstruction still requires them. Our midpoint-spline approach offers a more direct postprocessing route; it fits the Lanczos-derived cumulative midpoint data and differentiates the fit, without constructing a polynomial moment expansion or choosing a polynomial reference density and its scaling interval. 
% It still requires endpoint choices, and our comparisons with standard KPM--Jackson do not establish a performance advantage over spectrum-adaptive KPM.

\FloatBarrier
\subsection{Experimental Protocol}\label{sec:dos_experimental_protocol}

% Author checks before submission (not rendered): identify the exact driver
% that produced these figures. Record Lanczos reorthogonalization, probe
% normalization and whether methods share probes, KPM scaling bounds and
% Jackson convention, GP input scaling/hyperparameter fitting/jitter and
% whether the single GP used gp_monotonic. A previously supplied example
% used soft derivative constraints, so an unconstrained fit cannot be assumed.
% Also record the evaluation grid, padding/extension, negativity tolerance,
% software versions, seeds, and the provenance/interaction for Li7Nmax2.
% DOS_validation_eval.m returns absolute linfError and separately normalized
% normalizedLInfError; verify the plotting driver uses the relative
% max-error definition in eq:metrics_relative_linf, not either field alone.

\paragraph{Six DOS approximations.}
For each matrix, we compare the following six estimators, using the same labels as in the figures:
\begin{enumerate}
    \item \emph{Stochastic Lanczos}: the Gaussian-broadened quadrature DOS in~\eqref{eq:DOS_from_distr_reg}, averaged over random starting vectors.
    \item \emph{Ideal spline}: the DOS obtained by differentiating the monotone cubic CDOS interpolant constructed from ideal-starting-vector quadrature data, as motivated by~\eqref{eq:ideal_starting_vector}. This is an informative baseline that removes random-probe weighting from the underlying spectral measure. It still contains finite-$M$ quadrature and interpolation errors and is not the exact DOS. Its use of exact spectral information makes it a reference construction, rather than a practical competitor with the same information requirements.
    \item \emph{Midpoint spline}: our proposed DOS estimator~\eqref{eq:DOS_CS}, obtained by differentiating the monotone cubic interpolant of the averaged Lanczos midpoint data.
    \item \emph{Midpoint GP}: the derivative of a GP CDOS fit to the midpoint data using a single squared-exponential covariance kernel.
    \item \emph{Midpoint dual GP}: the analogous derivative using the sum of two squared-exponential kernels in~\eqref{eq:gp_dual_kernel}.
    \item \emph{KPM--Jackson}: the stochastic Chebyshev approximation with Jackson damping in~\eqref{eq:kpm_jackson_dos}.
\end{enumerate}
The GP DOS estimates are obtained by analytic differentiation of the fitted CDOS; equation~\eqref{eq:gp_dos} gives the Gaussian-likelihood case. In contrast to the monotone spline, the tested GP fits can yield negative DOS values; this distinction is retained in the evaluation rather than removed by clipping.

\paragraph{Common spectral resolution.}
For each matrix, we use its fixed Gaussian width from Table~\ref{tab:six_test_matrices} at every value of $M$. The four CDOS-derived DOS estimates and the KPM--Jackson estimate are convolved with this Gaussian before plotting and evaluation. The stochastic Lanczos curve is already the Gaussian convolution of its discrete quadrature measure, so it is not broadened a second time. All six curves are compared with the exact Gaussian DOS formed from the full eigenvalue set at that same width. Thus the comparisons follow the common-resolution convention in Section~\ref{sec:metrics}; the plotted spline and GP curves are the validation-broadened derivatives, not the raw derivatives. Their intrinsic interpolation or regression smoothing, and the Jackson damping of KPM, remain part of each method.

\paragraph{Iterations per probe, matrix--vector products, and repetitions.}
We vary
\begin{equation}\label{eq:dos_experiment_budgets}
    M\in\{9,15,30,60,90,120,150\},\qquad
    \nvec=R=5,
\end{equation}
where $R$ is the notation used for the number of random probes in the figures. The horizontal axes of the error-metric figures are labelled ``Iterations per probe, $M$''. For Lanczos-based methods, $M$ is the number of Lanczos steps per probe; for KPM it denotes the Chebyshev recurrence iteration count, corresponding to the polynomial degree, not a Lanczos iteration count. The common settings require a nominal total of $MR=5M$ matvecs per stochastic estimate, not $M$ matvecs in total. For example, $M=15$ and $M=90$ correspond to $75$ and $450$ matvecs, respectively, across the five probes in one repetition. Since $R$ is fixed, increasing $M$ increases the total matvec count proportionally. This comparison controls matrix access but does not equate elapsed time; orthogonalization, interpolation, GP fitting, and density evaluation incur different additional costs. The spectral information used by the ideal spline is not included in this practical cost comparison.

For each $(M,R)$ setting, we perform $B=10$ repetitions with different seeds for the random starting vectors. Five probes are averaged within each repetition to produce one DOS estimate; the ten repetitions are then used to summarize accuracy and variability. If $\widehat d_{\sigma,M}^{[b]}(t)$ denotes the compared DOS from repetition $b$ for a particular method, the comparison figures display
\begin{equation}\label{eq:dos_empirical_bands}
\begin{aligned}
    \overline d_{\sigma,M}(t)
    &=\frac1B\sum_{b=1}^{B}\widehat d_{\sigma,M}^{[b]}(t),\\
    s_d(t)
    &=\left[\frac{1}{B-1}\sum_{b=1}^{B}
       \left(\widehat d_{\sigma,M}^{[b]}(t)-\overline d_{\sigma,M}(t)\right)^2\right]^{1/2},
\end{aligned}
\end{equation}
with a shaded band $\overline d_{\sigma,M}(t)\pm s_d(t)$. These bands describe empirical variation between repetitions, not GP posterior credible intervals or confidence intervals for the mean. A band extending below zero does not by itself imply a negative estimated DOS; red curve segments specifically identify a negative approximation mean. Conversely, a nonnegative mean does not guarantee that every repetition was nonnegative.

For each error metric $E$, we first evaluate $E_b$ on each repetition and then plot its mean and sample standard deviation. In particular, the reported mean error is not the error of the mean DOS curve. The error figures use logarithmic vertical axes; as noted on the plots, error bars are omitted where the mean minus one standard deviation is nonpositive. Their absence therefore need not imply zero variability.

\paragraph{Undefined Jensen--Shannon values.}
The Jensen--Shannon divergence is evaluated only when the compared curve satisfies the probability-density requirements in Section~\ref{sec:metrics}. In particular, a GP repetition with materially negative DOS values has no valid divergence. The $k/10$ annotations report the number of finite, valid repetitions contributing to a plotted summary, with statistics computed only over that subset. When no repetition is valid, the plot marks the value as undefined and leaves a gap; it does not assign zero error. With only one valid repetition, a mean value is available but an empirical standard deviation is not. Since the plots use the general criterion of finiteness, a count annotation on another method should also be read as an availability count, not automatically as evidence of negativity. Comparisons involving different validity counts must be interpreted cautiously as a small conditional mean divergence does not compensate for invalid outputs in the remaining repetitions.

\paragraph{Presentation of results.}
Every DOS comparison uses the same $2\times3$ panel order: stochastic Lanczos, ideal spline, and midpoint spline in the top row; midpoint GP, midpoint dual GP, and KPM--Jackson in the bottom row. The black curve is the exact Gaussian DOS. Panel axis limits can differ, so the size of an apparent discrepancy should be judged from its numerical scale, not just its height on the page. Error summaries use a $2\times2$ arrangement: relative $L^\infty$ and relative $L^2$ errors above, and Jensen--Shannon divergence and cosine error below. We discuss the matrices in the order used in Figures~\ref{fig:six_matrix_eigenvalue_spectra} and~\ref{fig:six_matrix_exact_DOS}, moving from interior oscillations to fine-feature resolution and then to approximation with few matvecs of strongly clustered spectra.

\FloatBarrier
\subsection{Benzene: Recovering the Interior Density}\label{sec:results_benzene}

\begin{figure}[tbp]
\centering
\includegraphics[width=\textwidth]{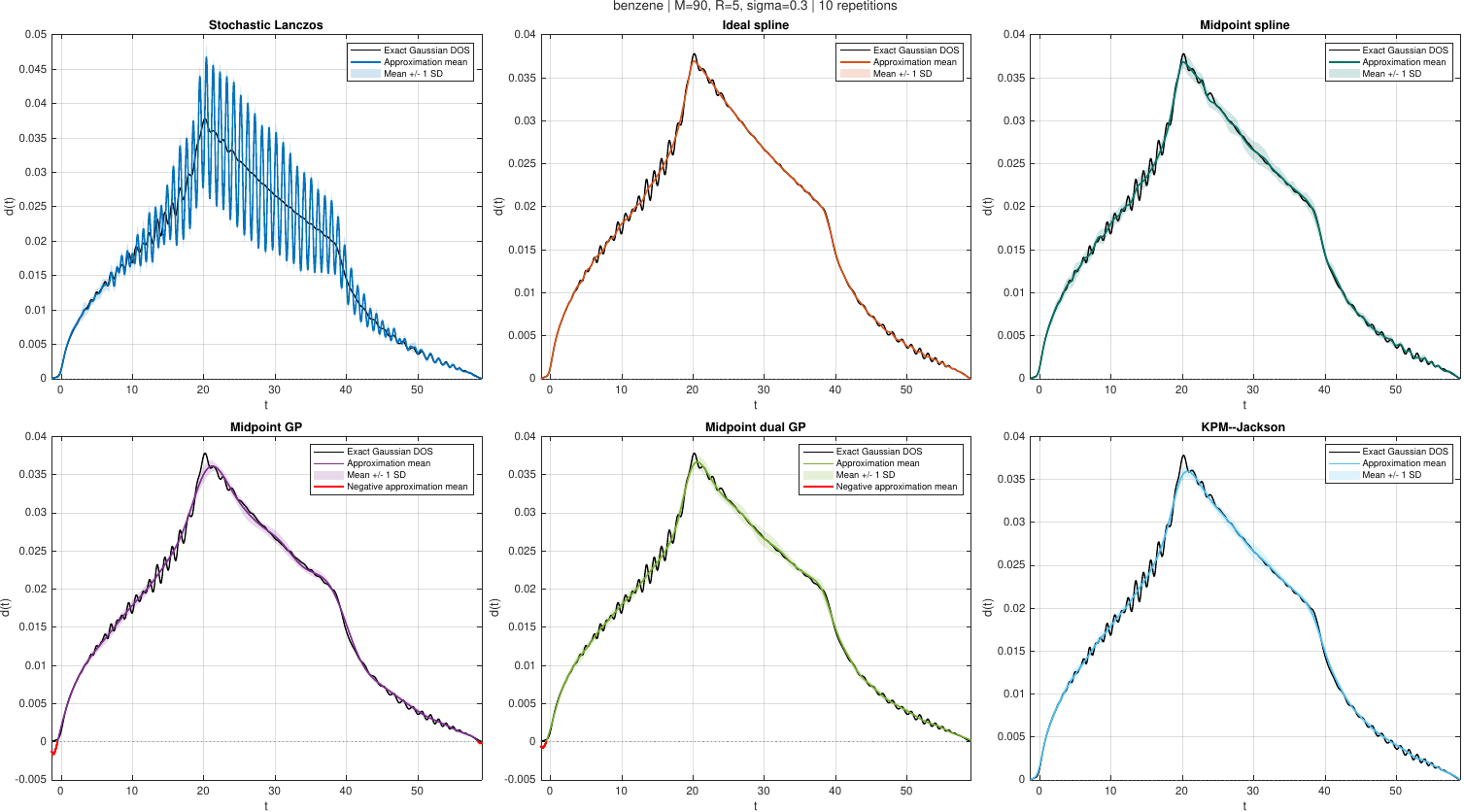}
\caption{\texttt{benzene}: DOS comparisons at $M=90$. The panel order and mean $\pm$ one empirical standard-deviation bands follow Section~\ref{sec:dos_experimental_protocol}. Stochastic Lanczos exhibits interior oscillations, whereas the midpoint spline follows the overall profile without them. Red segments identify a negative approximation mean.}
\label{fig:benzene_dos_m90}
\end{figure}

Figure~\ref{fig:benzene_dos_m90} shows a clear difference between representing a density by broadened quadrature nodes and differentiating a fitted cumulative distribution. At $M=90$, stochastic Lanczos captures the low- and high-energy portions of the spectrum well, but its interior contains large, closely spaced oscillations about the reference curve. The midpoint spline follows the dominant peak, the descending shoulder, and the overall tail without these quadrature-node spikes. Both GP fits also follow the broad profile closely, although the single-kernel GP rounds the dominant peak more strongly. KPM--Jackson gives a good broad profile but also smooths some local detail. The ideal spline provides a particularly accurate reference curve at this value of $M$.

This Lanczos behavior is consistent with the typically faster convergence of extremal Ritz values and with insufficient overlap between Gaussian kernels centered on the interior quadrature nodes~\cite{Saad2011,LSY}. The relevant gaps are those between the \emph{Ritz values}, not changes in the true eigenvalue spacing, which is fixed. Increasing $M$ refines the quadrature representation and can reduce these artificial oscillations at a fixed $\sigma$. Accurate DOS estimation nevertheless concerns the weighted spectral measure, not the convergence of every Ritz value to an individual eigenvalue.

\begin{figure}[tbp]
\centering
\includegraphics[width=\textwidth]{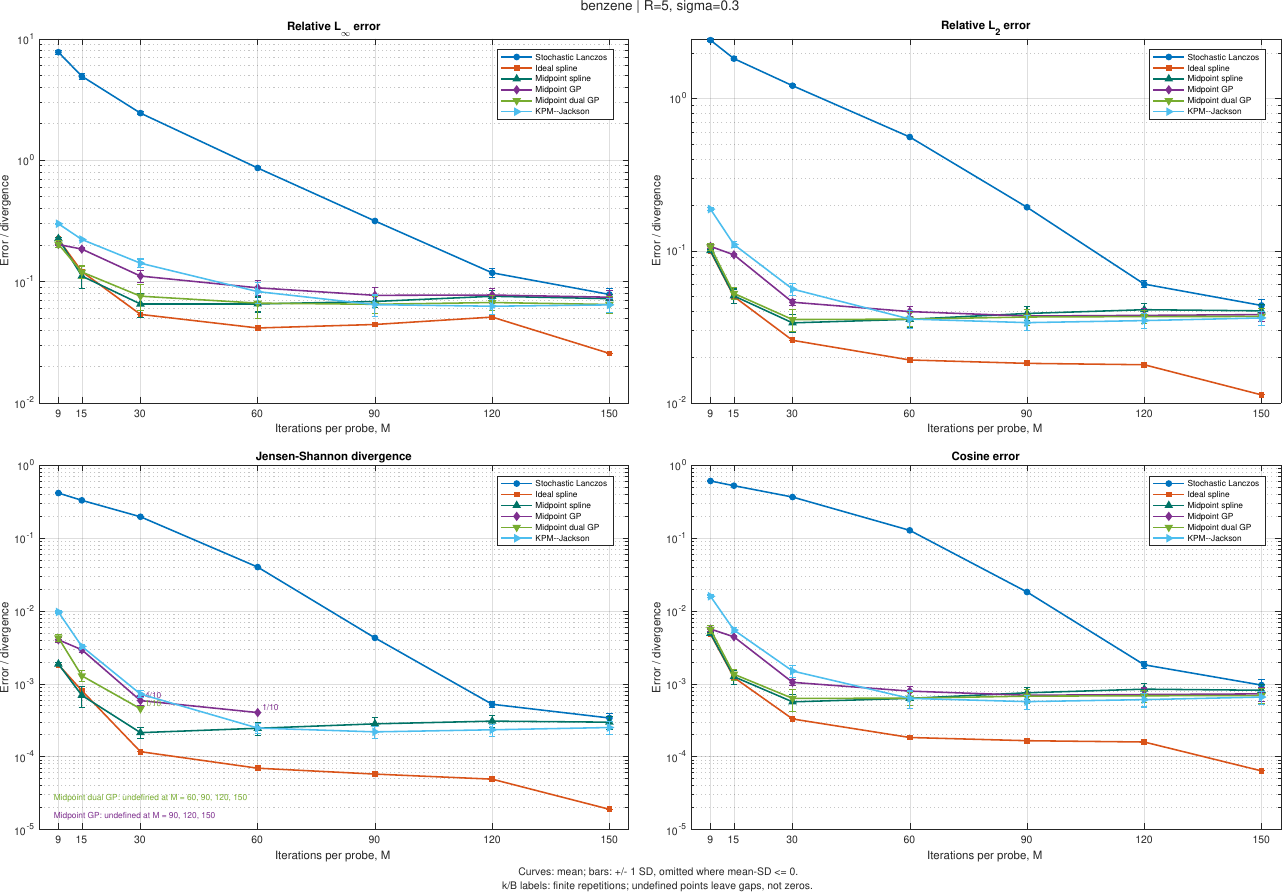}
\caption{\texttt{benzene}: error means and empirical standard deviations over ten repetitions as the number of iterations per probe $M$ varies, with $R=5$ and $\sigma=0.3$. The top row shows relative $L^\infty$ and $L^2$ errors; the bottom row shows Jensen--Shannon divergence and cosine error. Finite-repetition counts and undefined values are interpreted as in Section~\ref{sec:dos_experimental_protocol}.}
\label{fig:benzene_dos_errors}
\end{figure}

The error curves in Figure~\ref{fig:benzene_dos_errors} confirm the substantial improvement of stochastic Lanczos with increasing $M$. The midpoint spline reaches much lower errors with a moderate number of matvecs, while the practical spline, both GP variants, and KPM curves become close in the norm and cosine measures at the largest values of $M$. Their errors need not decrease strictly at every step. The ideal spline generally remains more accurate, showing that increasing $M$ with only five random probes does not remove all of the discrepancy from the ideal construction.

The GP fits' good visual agreement also requires a qualification: small negative portions of their means are marked in red near the spectral endpoints of Figure~\ref{fig:benzene_dos_m90}. The single GP's Jensen--Shannon divergence is undefined in all ten repetitions at $M=90,120,150$, while the dual GP's is undefined at $M=60,90,120,150$. The midpoint spline therefore combines competitive shape accuracy with the nonnegativity needed for a density interpretation, rather than achieving only a close curve fit.

\FloatBarrier
\subsection{HubbardL8: The Effect of Increasing the Iterations per Probe}\label{sec:results_hubbard}

\begin{figure}[tbp]
\centering
\includegraphics[width=\textwidth]{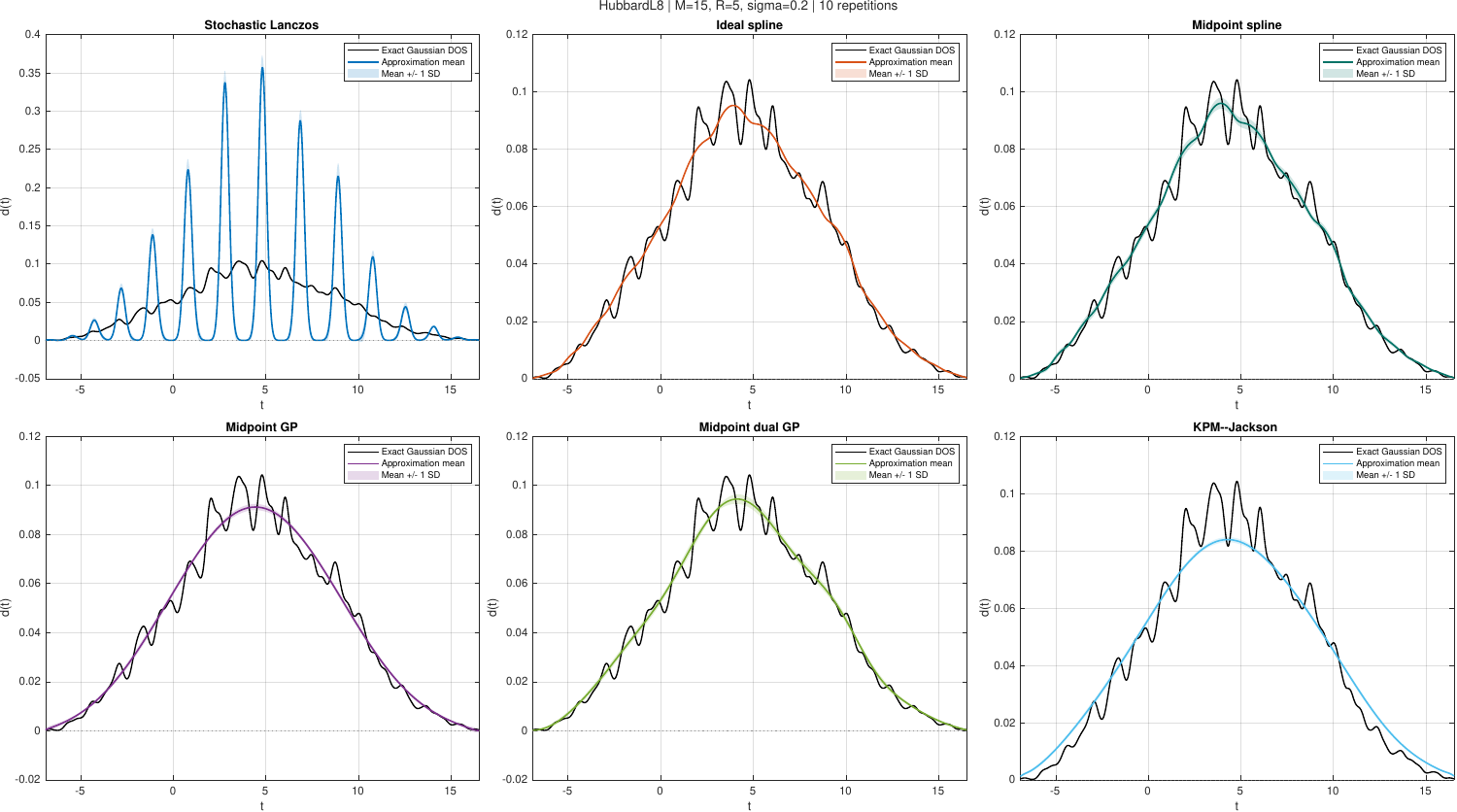}
\caption{\texttt{HubbardL8}: DOS comparisons at $M=15$. The midpoint spline recovers the broad envelope despite using few matvecs, while stochastic Lanczos produces separated, excessively high peaks. The Lanczos panel uses a different vertical scale. Bands show mean $\pm$ one empirical standard deviation over ten repetitions.}
\label{fig:hubbard_dos_m15}
\end{figure}

The Hubbard example demonstrates both the advantage of cumulative fitting when few matvecs are available and a regime in which stochastic Lanczos becomes highly competitive. At $M=15$ in Figure~\ref{fig:hubbard_dos_m15}, the Lanczos estimate is a sequence of separated peaks with amplitudes much larger than those of the reference DOS. Both spline constructions recover the broad central envelope, though neither resolves all of the genuine local maxima with so few matvecs. The GP fits likewise smooth over much of the local structure but give competitive integrated accuracy; in Table~\ref{tab:dos-relative-l2}, the dual GP has the smallest practical-method mean relative $L^2$ error for this setting. KPM--Jackson also produces a smooth envelope, with noticeably lowered central peak heights.

\begin{figure}[tbp]
\centering
\includegraphics[width=\textwidth]{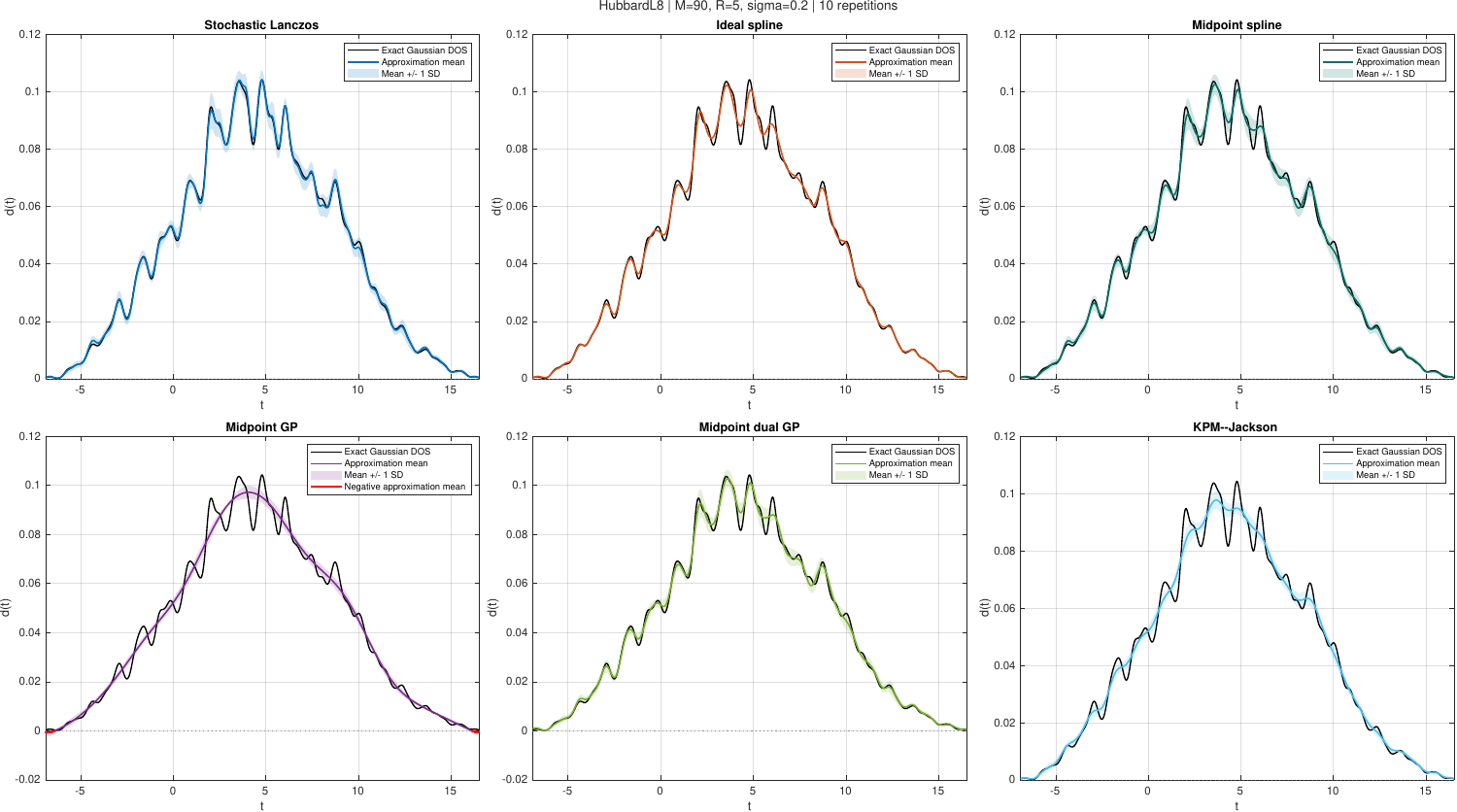}
\caption{\texttt{HubbardL8}: the corresponding comparison at $M=90$. Stochastic Lanczos now resolves the main local peaks and valleys. The midpoint spline and dual GP also recover substantial local detail, whereas the single-kernel GP and KPM--Jackson produce smoother profiles.}
\label{fig:hubbard_dos_m90}
\end{figure}

At $M=90$ in Figure~\ref{fig:hubbard_dos_m90}, stochastic Lanczos resolves the principal local peaks and valleys and is among the most accurate practical estimators. The midpoint spline also captures these features well and remains close to the ideal spline, while the dual GP gives a comparably detailed fit. KPM--Jackson retains more smoothing, and the single-kernel GP largely reproduces only the broad envelope. Thus the spline's advantage is not that stochastic Lanczos cannot recover the DOS, but that useful shape information is obtained before enough quadrature nodes are available for a well-resolved Gaussian sum.

\begin{figure}[tbp]
\centering
\includegraphics[width=\textwidth]{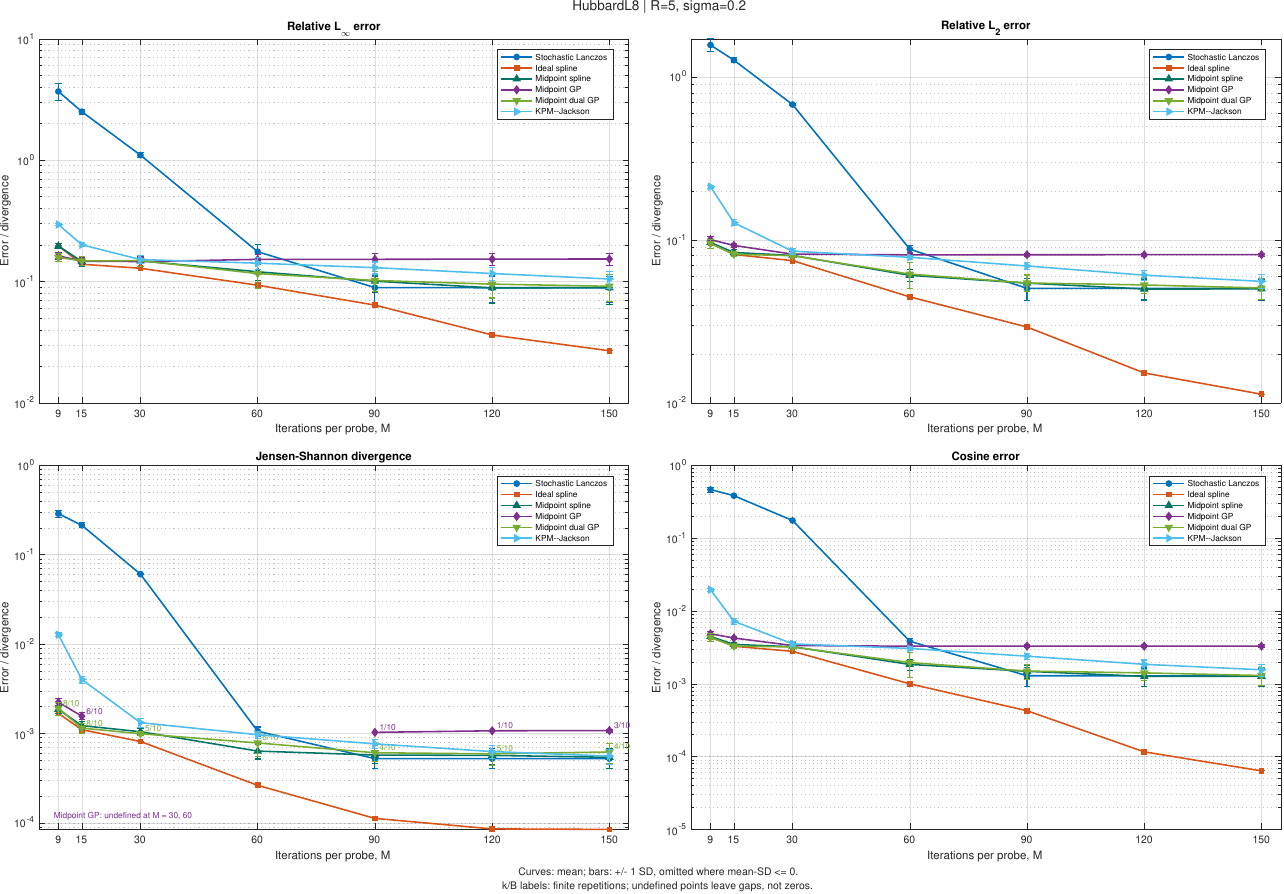}
\caption{\texttt{HubbardL8}: error means and empirical standard deviations versus the number of iterations per probe $M$, with $R=5$ and $\sigma=0.2$. The panel order is as in Figure~\ref{fig:benzene_dos_errors}. The improvement of stochastic Lanczos brings it close to the midpoint spline and dual GP with more matvecs.}
\label{fig:hubbard_dos_errors}
\end{figure}

Figure~\ref{fig:hubbard_dos_errors} supports this transition. The Lanczos errors drop rapidly as $M$ increases, and by about $M=90$ they are close to those of the midpoint spline and dual GP. Small differences between these practical estimators depend on the metric and occur alongside repetition variability, so the figures do not establish one strict ranking for every measure. The ideal spline continues to improve with $M$, whereas the single-kernel GP errors largely level off. Its Jensen--Shannon divergence is undefined in all ten repetitions at $M=30$ and $60$, and several other plotted values use only the valid subset of repetitions. The useful comparison is therefore both accuracy with few matvecs and the accuracy attained with more matvecs, not a claim of uniform dominance by one estimator.

\FloatBarrier
\subsection{HeisenbergL16: More Iterations per Probe Can Still Be Insufficient}\label{sec:results_heisenberg}

\begin{figure}[tbp]
\centering
\includegraphics[width=\textwidth]{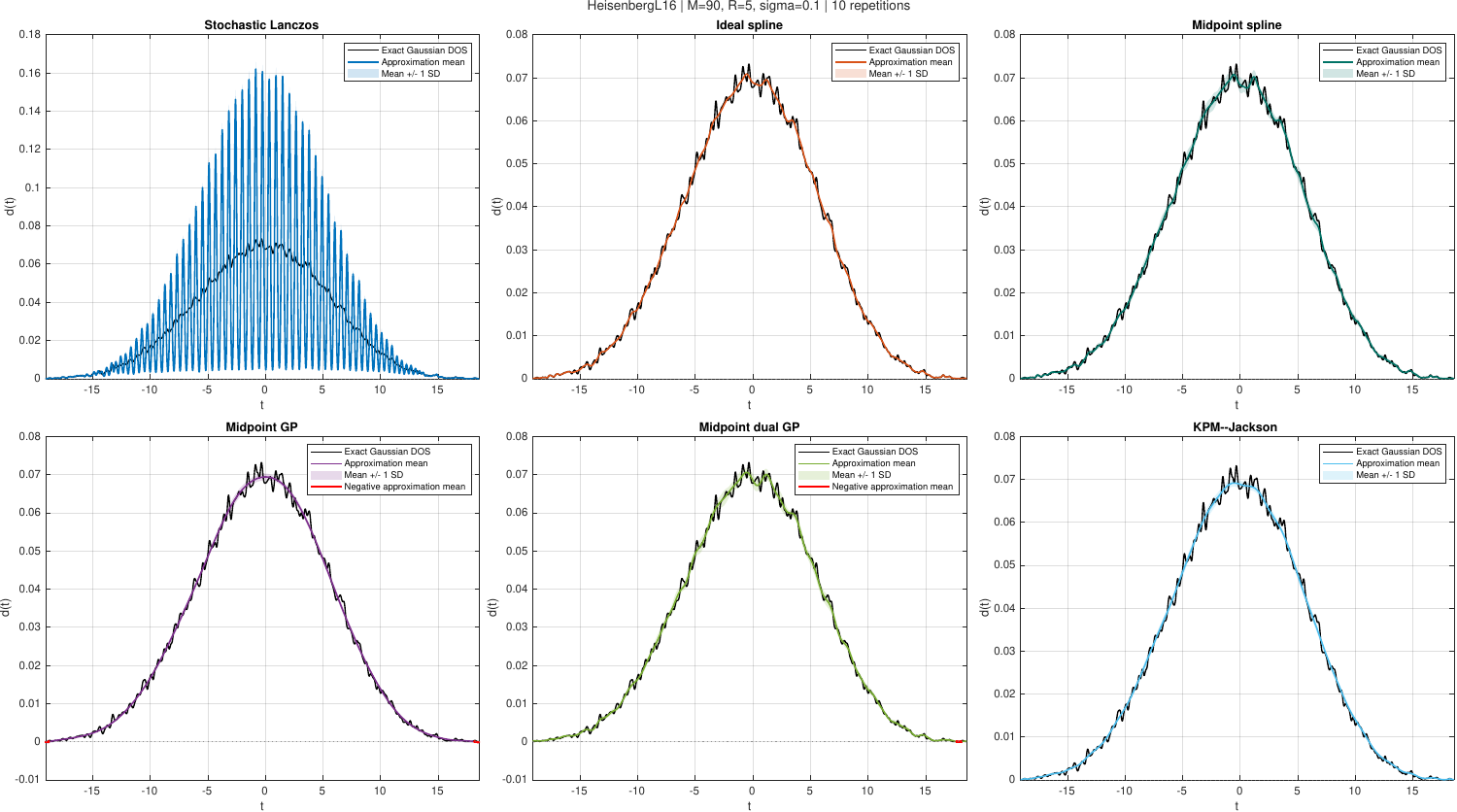}
\caption{\texttt{HeisenbergL16}: DOS comparisons at $M=90$. The Lanczos estimate retains large interior oscillations; note its larger vertical scale. The midpoint spline follows the reference envelope and local variations closely. Bands show empirical variation over ten repetitions, and red segments mark a negative mean.}
\label{fig:heisenberg_dos_m90}
\end{figure}

The same $M=90$ iterations per probe does not give an equally resolved Lanczos DOS for every problem. Figure~\ref{fig:heisenberg_dos_m90} shows pronounced oscillations throughout the central region of the Heisenberg spectrum at $\sigma=0.1$. The peaks in the Lanczos panel substantially exceed the reference height, while its tail regions are better represented. In contrast, the midpoint spline closely follows the broad central density and several of its local variations, with a profile similar to that of the ideal spline. Both GP fits recover the broad density accurately, with the dual GP retaining somewhat more local variation and the single GP producing a smoother profile. Both also have red-marked negative segments in the tails. KPM--Jackson similarly follows the broad envelope while smoothing finer variations.

Compared with the Hubbard example, this result shows that the number of matvecs required by Lanczos depends on the spectral distribution and the chosen Gaussian resolution, not on $M$ or matrix dimension alone. At the displayed value of $M$, broadening the discrete quadrature measure still exposes its individual interior nodes, whereas interpolating cumulative information gives a useful density profile. Table~\ref{tab:dos-relative-l2} quantifies this difference: at $M=90$, the mean relative $L^2$ error is $0.7600$ for stochastic Lanczos and $0.03563$ for the midpoint spline, a factor of approximately $21.3$. The single GP and KPM--Jackson have slightly smaller mean $L^2$ errors than the midpoint spline here, so the visual recovery of local structure should not be interpreted as a strict ranking under this integrated metric.

% An error-versus-M figure for HeisenbergL16 was not supplied with this set.

\FloatBarrier
\subsection{Li7Nmax2: Preserving Resolved Local Structure}\label{sec:results_li7}

\begin{figure}[tbp]
\centering
\includegraphics[width=\textwidth]{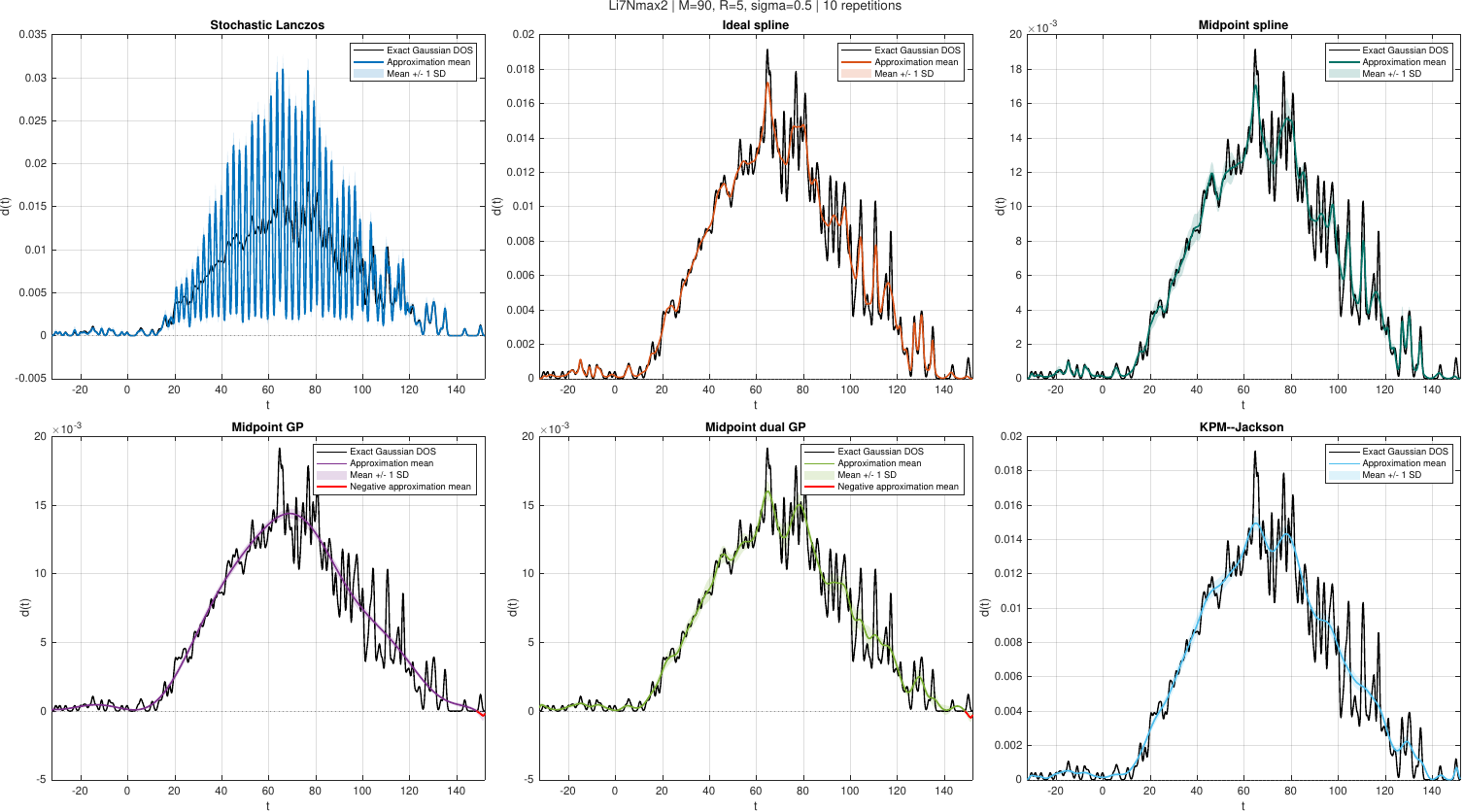}
\caption{\texttt{Li7Nmax2}: DOS comparisons at $M=90$. The midpoint spline preserves more of the reference peak structure than the GP and KPM--Jackson fits, although some narrow features remain unresolved. The stochastic Lanczos panel has a larger vertical scale. Bands show empirical variation over ten repetitions.}
\label{fig:li7_dos_m90}
\end{figure}

For \texttt{Li7Nmax2}, the challenge is not only to recover a broad envelope but also to retain the many local features that remain visible at $\sigma=0.5$. Figure~\ref{fig:li7_dos_m90} shows that stochastic Lanczos at $M=90$ again produces large interior oscillations unrelated to the reference peak pattern. The midpoint spline instead follows several of the prominent local peaks and shoulders and is close to the ideal spline. The dual GP and KPM--Jackson provide reasonable large-scale profiles but merge or attenuate more of the local peaks. The single-kernel GP is smoother still and mainly captures the envelope.

\begin{figure}[tbp]
\centering
\includegraphics[width=\textwidth]{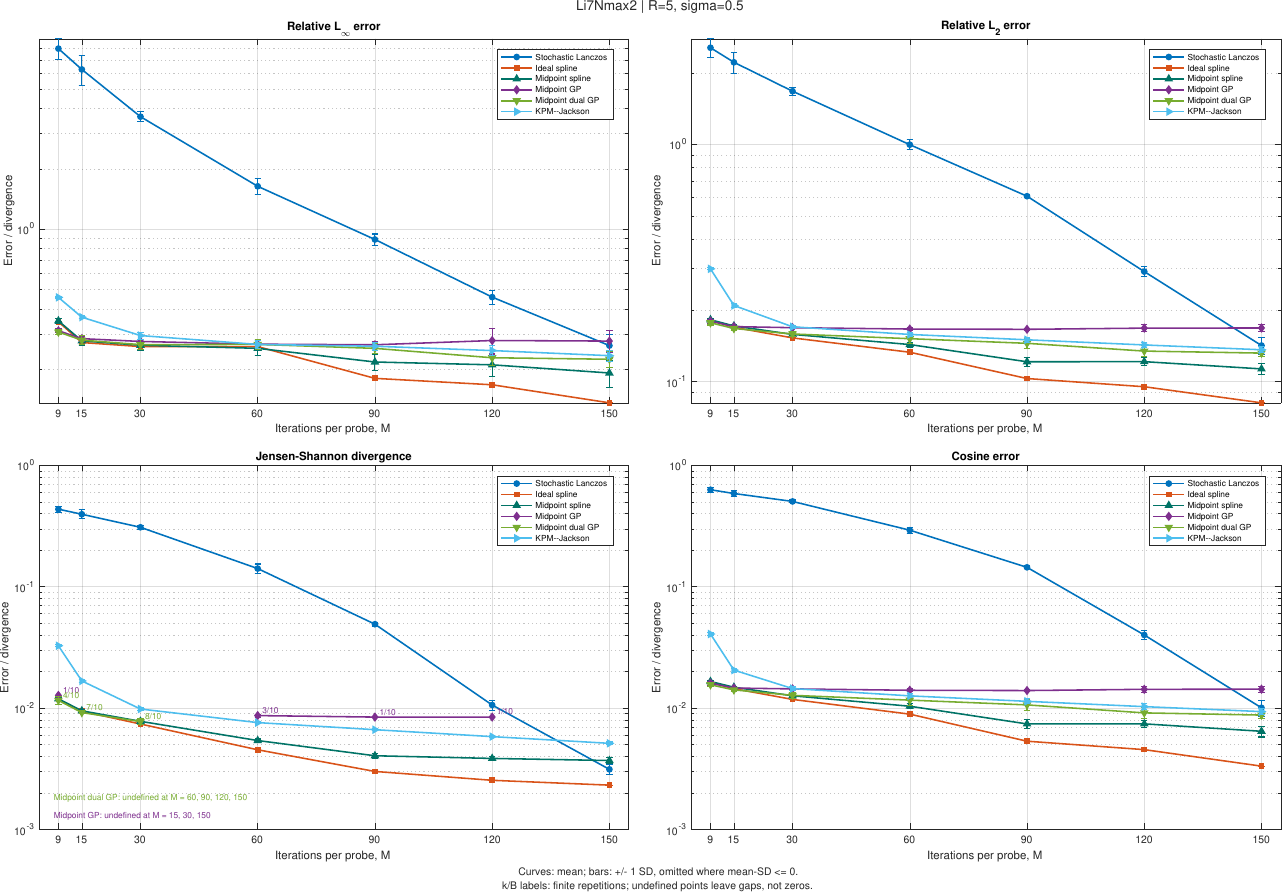}
\caption{\texttt{Li7Nmax2}: error means and empirical standard deviations versus the number of iterations per probe $M$, with $R=5$ and $\sigma=0.5$. The panel order is as in Figure~\ref{fig:benzene_dos_errors}. Undefined GP divergences are shown as gaps; a finite-subset mean must be read together with its validity count.}
\label{fig:li7_dos_errors}
\end{figure}

This is a resolution advantage of the midpoint spline, not evidence that all fine features are recovered as some narrow peaks remain underestimated or merged by both spline variants. Figure~\ref{fig:li7_dos_errors} places the midpoint spline among the best practical estimators over much of the tested range and close to the ideal reference. The single GP is also competitive in the relative norm measures at small $M$ and then changes comparatively little as $M$ increases. Its Jensen--Shannon divergence is undefined in all ten repetitions at $M=15,30,150$ and is computed from only a subset of repetitions at some other values of $M$. The dual GP has no valid Jensen--Shannon values at $M=60,90,120,150$.

Stochastic Lanczos improves substantially as $M$ increases. At $M=150$, its plotted Jensen--Shannon divergence is slightly below that of the midpoint spline, even though its relative norm and cosine errors remain higher. This metric dependence reinforces the need to consider the density plots together with several error measures. The spline's main benefit here is a favorable balance of local resolution, overall agreement, and density validity with a moderate number of matvecs.

\FloatBarrier
\subsection{nd3k: Accuracy with Few Matvecs and Complementary Metrics}\label{sec:results_nd3k}

\begin{figure}[tbp]
\centering
\includegraphics[width=\textwidth]{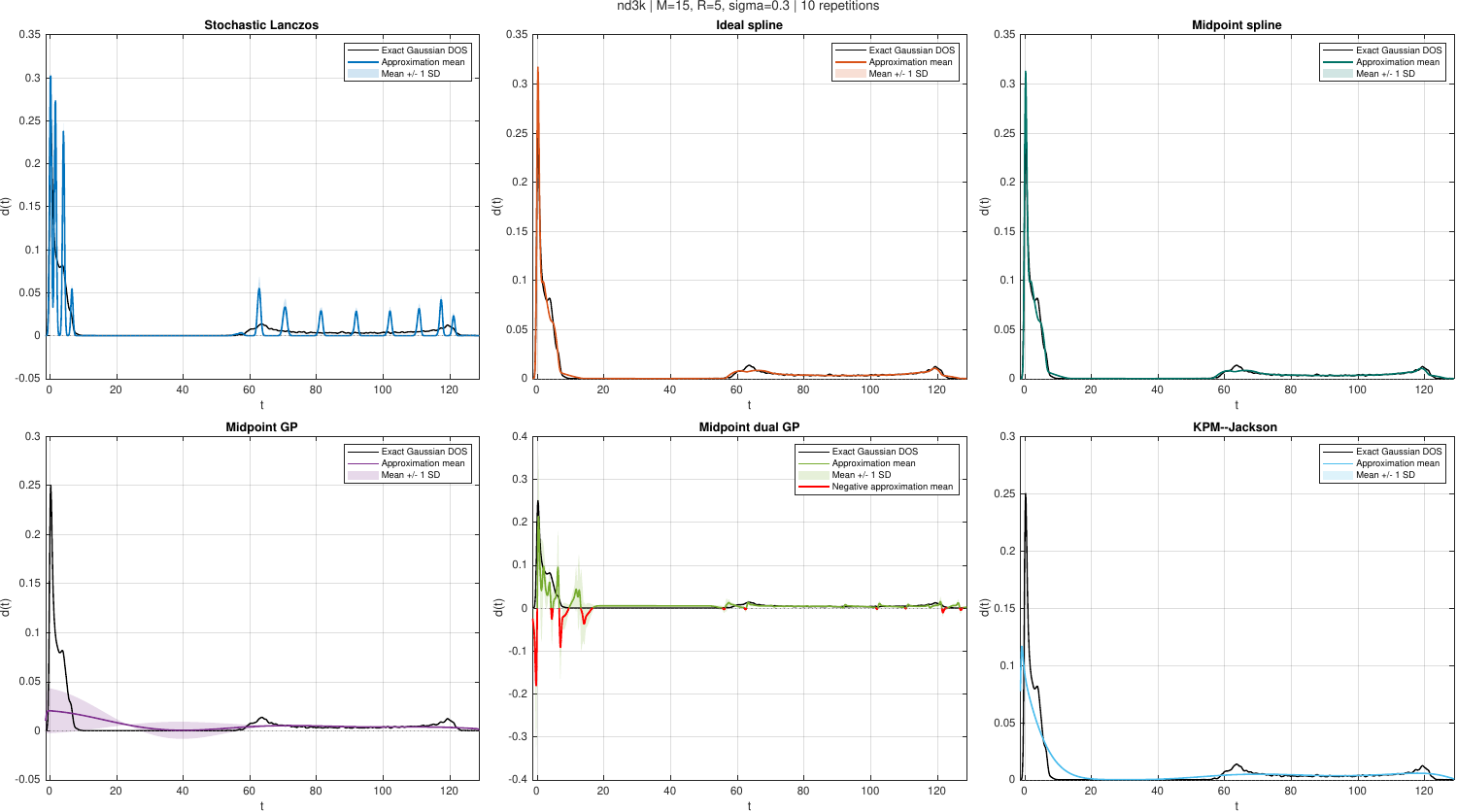}
\caption{\texttt{nd3k}: DOS comparisons at $M=15$. The midpoint spline captures the main spectral groups with few matvecs. The GP panels show negative mean regions, and the dual GP has substantial variability. Vertical scales differ between panels.}
\label{fig:nd3k_dos_m15}
\end{figure}

The \texttt{nd3k} matrix highlights the usefulness of cumulative interpolation when the spectrum is strongly clustered. At only $M=15$, Figure~\ref{fig:nd3k_dos_m15} shows that the midpoint spline already identifies the dominant low-energy concentration, the large gap, and the smaller high-energy density features. It follows the ideal spline closely, although the dominant peak height is not exact. Stochastic Lanczos exhibits separated peaks within the spectral groups, while KPM--Jackson strongly attenuates the dominant peak. Both GP variants have negative mean regions; the dual GP additionally shows substantial repeating variability.

\begin{figure}[tbp]
\centering
\includegraphics[width=\textwidth]{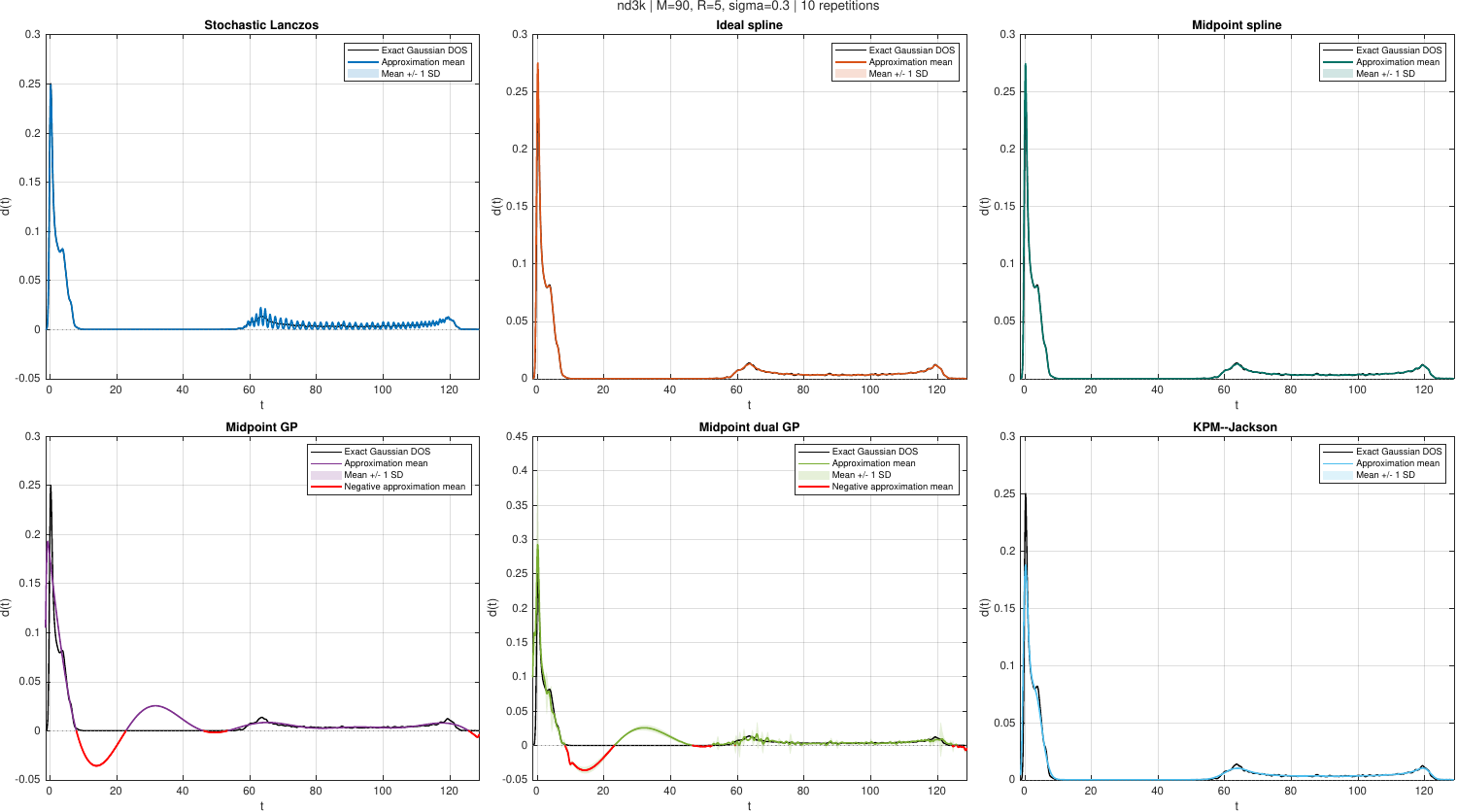}
\caption{\texttt{nd3k}: the corresponding comparison at $M=90$. Stochastic Lanczos captures the dominant low-end peak accurately but retains oscillations in the higher-energy group. The midpoint spline gives a smoother overall fit with a residual error in the largest peak height.}
\label{fig:nd3k_dos_m90}
\end{figure}

At $M=90$ in Figure~\ref{fig:nd3k_dos_m90}, the midpoint spline gives a close overall fit, with its most visible discrepancy near the height of the largest peak. Stochastic Lanczos matches this low-end peak especially well but retains oscillations in the much smaller high-energy density. The GP estimates still exhibit excursions below zero and artificial structure in the gap. Thus adding quadrature information does not by itself make a GP fit a valid density.

\begin{figure}[tbp]
\centering
\includegraphics[width=\textwidth]{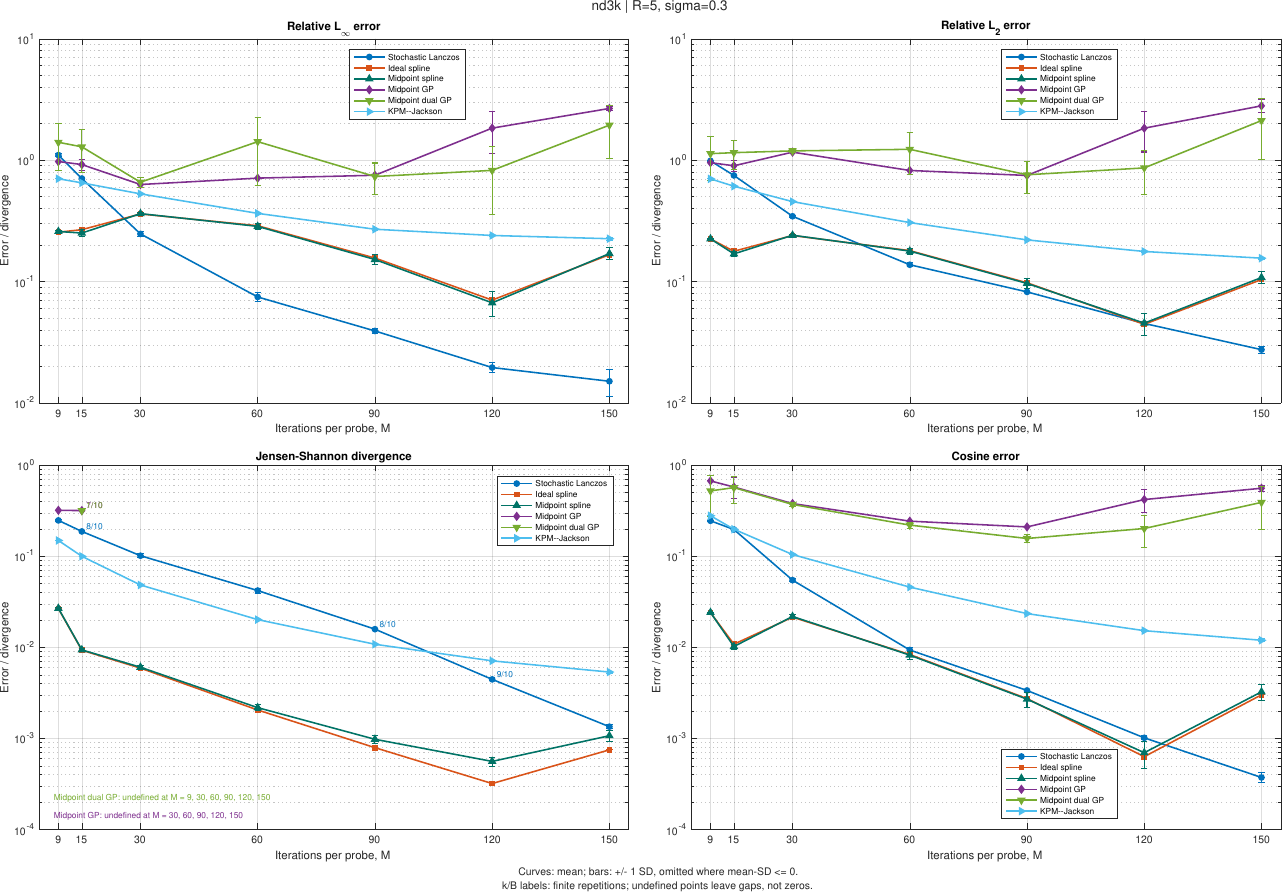}
\caption{\texttt{nd3k}: error means and empirical standard deviations versus the number of iterations per probe $M$, with $R=5$ and $\sigma=0.3$. The panel order is as in Figure~\ref{fig:benzene_dos_errors}. Norm errors and Jensen--Shannon divergence can favor different estimators. The spline errors are not monotone in $M$, and GP divergence availability is severely limited.}
\label{fig:nd3k_dos_errors}
\end{figure}

Figure~\ref{fig:nd3k_dos_errors} illustrates why the preferred estimator depends on what is measured. At $M=90$, stochastic Lanczos has a lower relative maximum error than the midpoint spline, reflecting its accurate representation of the dominant peak. This is not a failure of the $L^\infty$ metric: it correctly emphasizes the largest absolute discrepancy. However, that normalization makes smaller-amplitude oscillations elsewhere less influential. The Jensen--Shannon divergence, which compares normalized mass distributions, favors the spline more strongly at this value of $M$ and reveals the benefit of its smoother representation of the high-energy group. Cosine error is $L^2$-weighted and therefore need not distinguish those smaller features as strongly.

The spline and ideal-spline error curves track one another closely on this problem, including their nonmonotone behavior. In particular, several errors increase from $M=120$ to $150$. Stochastic Lanczos, by contrast, continues to improve and achieves lower norm and cosine errors at the largest tested value of $M$, while the spline retains a smaller Jensen--Shannon divergence. These observations support an advantage with few matvecs and for probability-shape agreement, rather than monotone convergence or uniform superiority in every metric. The validity annotations also matter: the single GP has no finite Jensen--Shannon values at $M=30,60,90,120,150$, and the dual GP has none at $M=9,30,60,90,120,150$. At the remaining iteration counts, the displayed GP divergence statistics can still be based on fewer than ten valid repetitions.

\FloatBarrier
\subsection{fv1: Capturing an Upper-End Concentration at Small Cost}\label{sec:results_fv1}

\begin{figure}[tbp]
\centering
\includegraphics[width=\textwidth]{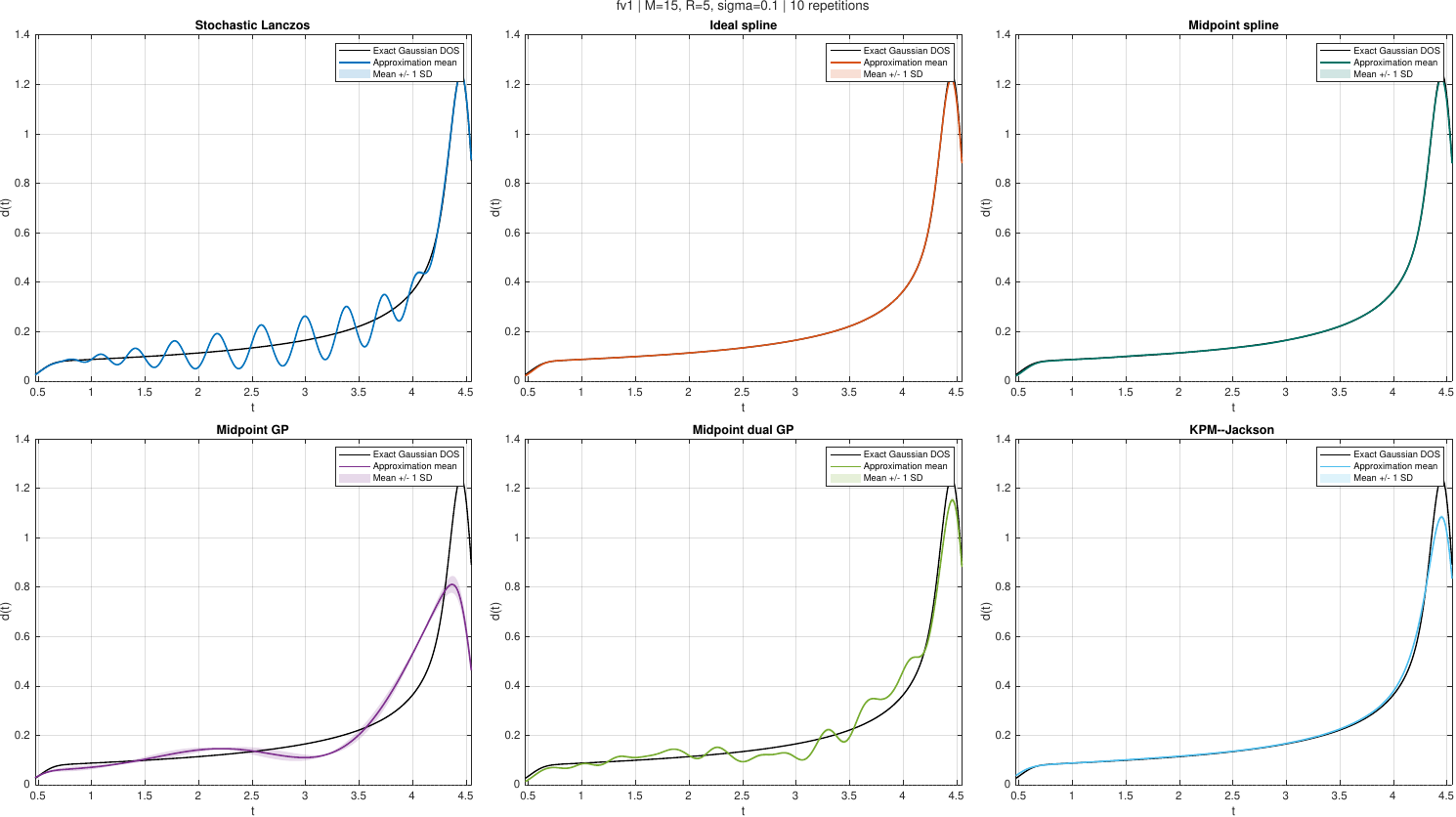}
\caption{\texttt{fv1}: DOS comparisons at $M=15$. The midpoint spline follows the smooth background and the dominant upper-end peak closely. Stochastic Lanczos retains interior oscillations, while KPM--Jackson and both GP variants underestimate the peak. Bands show empirical variation over ten repetitions.}
\label{fig:fv1_dos_m15}
\end{figure}

The final example demonstrates accuracy with few matvecs for a density concentrated near the opposite end of the spectrum. In Figure~\ref{fig:fv1_dos_m15}, with $M=15$ and $\sigma=0.1$, the midpoint spline is already visually close to the ideal spline and the exact Gaussian DOS. It captures the smooth background, the rise toward the upper spectral edge, and both the location and height of the dominant peak.

Stochastic Lanczos also reproduces the upper-end peak well but exhibits visible oscillations across the lower-density interior. KPM--Jackson gives a smooth and otherwise close profile, with the largest difference in the attenuated peak height. The dual GP produces artificial undulations along the rising density and underestimates the peak, while the single GP shifts and substantially lowers the peak and distorts the background. With so few matvecs, the midpoint spline therefore offers a particularly effective combination of overall shape and peak accuracy.

The mean relative $L^2$ errors in Table~\ref{tab:dos-relative-l2} support this advantage with few matvecs: at $M=15$, the midpoint spline gives $0.01461$, compared with $0.1356$ for stochastic Lanczos and $0.08711$ for KPM--Jackson. Increasing the number of iterations per probe to $M=90$ brings the latter two errors down to $0.01244$ and $0.01246$, respectively, slightly below the spline's $0.01342$. Thus the spline already delivers at $M=15$ an accuracy close to that attained by these alternatives at six times the number of iterations per probe (and hence six times the total matvec count at fixed $R$).

% An error-versus-M figure for fv1 was not supplied with this set.

\FloatBarrier
\subsection{Quantitative Summary and Implications}\label{sec:dos_results_summary}

Tables~\ref{tab:dos-relative-l2} and~\ref{tab:dos-spline-rank-gap} provide complementary summaries of the mean relative $L^2$ errors at $M=15$ and $M=90$ iterations per probe, with $R=5$ probes. The first gives the reported error values for all six matrices and all methods; the second shows the midpoint spline's rank and its error gap relative to the best practical estimator in each setting. Together, they distinguish achieving the smallest error from remaining consistently close to it. The ideal spline is retained as a reference but excluded from practical-method rankings. Both tables summarize means over ten repetitions, so small differences should be interpreted alongside the empirical variability where available, not as statistically established superiority.

\begin{table}[!htbp]
\centering
\small
\setlength{\tabcolsep}{4pt}
\caption{Mean relative $L^2$ errors over ten repetitions with $R=5$ probes and the matrix-specific Gaussian widths in Table~\ref{tab:six_test_matrices}. The two settings use $M=15$ and $M=90$ iterations per probe, corresponding to nominal totals of $75$ and $450$ matvecs per repetition for the practical estimators. Values are relative errors, not percentages. Bold entries identify the smallest practical-method mean in each row, excluding the ideal-spline reference; they do not indicate statistical significance.}
\label{tab:dos-relative-l2}
\begin{tabular}{@{}lr|rrrrr@{}}
    \toprule
    Matrix
    & \shortstack{Ideal spline\\(reference)}
    & \shortstack{Stochastic\\Lanczos}
    & \shortstack{Midpoint\\spline}
    & \shortstack{Midpoint\\GP}
    & \shortstack{Midpoint\\dual GP}
    & \shortstack{KPM--\\Jackson} \\
    \midrule
    \multicolumn{7}{l}{$M=15$} \\
    \addlinespace[2pt]
    \texttt{benzene}       & 0.05004  & 1.835   & \textbf{0.05049} & 0.09404 & 0.05243 & 0.1100  \\
    \texttt{HubbardL8}     & 0.08141  & 1.270   & 0.08375 & 0.09253 & \textbf{0.08171} & 0.1276  \\
    \texttt{HeisenbergL16} & 0.04403  & 2.265   & 0.04524 & \textbf{0.03905} & 0.03921 & 0.1305  \\
    \texttt{Li7Nmax2}      & 0.1697   & 2.228   & 0.1720  & 0.1708  & \textbf{0.1677} & 0.2098  \\
    \texttt{nd3k}          & 0.1774   & 0.7476  & \textbf{0.1689} & 0.8984  & 1.152   & 0.6096  \\
    \texttt{fv1}           & 0.01183  & 0.1356  & \textbf{0.01461} & 0.3384 & 0.1225  & 0.08711 \\
    \midrule
    \multicolumn{7}{l}{$M=90$} \\
    \addlinespace[2pt]
    \texttt{benzene}       & 0.01818  & 0.1933  & 0.03870 & 0.03722 & 0.03668 & \textbf{0.03369} \\
    \texttt{HubbardL8}     & 0.02921  & \textbf{0.05043} & 0.05430 & 0.08115 & 0.05445 & 0.06921 \\
    \texttt{HeisenbergL16} & 0.02455  & 0.7600  & 0.03563 & 0.03375 & 0.03572 & \textbf{0.03366} \\
    \texttt{Li7Nmax2}      & 0.1032   & 0.6065  & \textbf{0.1214} & 0.1664 & 0.1452  & 0.1503  \\
    \texttt{nd3k}          & 0.09755  & \textbf{0.08208} & 0.09638 & 0.7476 & 0.7572  & 0.2200  \\
    \texttt{fv1}           & 0.005850 & \textbf{0.01244} & 0.01342 & 0.1552 & 0.08214 & 0.01246 \\
    \bottomrule
\end{tabular}
\end{table}

\paragraph{Accuracy with few iterations per probe.}
At $M=15$, Table~\ref{tab:dos-relative-l2} shows that the midpoint spline has a lower mean relative $L^2$ error than both stochastic Lanczos and KPM--Jackson for every matrix. Relative to stochastic Lanczos, its errors are smaller by factors ranging from approximately $4.4$ for \texttt{nd3k} to $50.1$ for \texttt{HeisenbergL16}. The midpoint spline has the smallest mean for \texttt{benzene}, \texttt{nd3k}, and \texttt{fv1} among all practical methods; dual GP has the smallest mean for \texttt{HubbardL8} and \texttt{Li7Nmax2}, and single GP has the smallest mean for \texttt{HeisenbergL16}. Thus the spline's advantage at this iteration count is clear relative to the direct Lanczos and polynomial estimators, although GP regression can be competitive or more accurate under this metric.

\paragraph{The effect of increasing the iteration count.}
Increasing to $M=90$ uses six times as many matvecs at the fixed probe count and changes the ordering of the methods. The midpoint spline has the smallest practical-method mean only for \texttt{Li7Nmax2}, which is not the largest matrix but has more local detail features in its DOS as shown in Figure~\ref{fig:li7_dos_m90}. Stochastic Lanczos has the smallest mean for \texttt{HubbardL8}, \texttt{nd3k}, and \texttt{fv1}, while KPM--Jackson has the smallest mean for \texttt{benzene} and \texttt{HeisenbergL16}. The single GP is nevertheless close to the KPM--Jackson minimum on \texttt{HeisenbergL16}, with means of $0.03375$ and $0.03366$, respectively. The spline's mean errors remain approximately five times smaller than those of stochastic Lanczos for \texttt{benzene} and \texttt{Li7Nmax2}, and more than twenty times smaller for \texttt{HeisenbergL16}. The benefit of additional iterations is therefore strongly problem-dependent. Moreover, some differences between methods are very small: on \texttt{fv1}, the Lanczos and KPM--Jackson means are $0.01244$ and $0.01246$. This motivates examining error gaps as well as first-place results.

\paragraph{Consistency relative to the best practical estimator.}
To quantify that distinction, Table~\ref{tab:dos-spline-rank-gap} pairs the midpoint spline's rank among the five practical estimators with its percentage excess over the smallest mean relative $L^2$ error. For each matrix and each listed value of $M$, let $\overline E_{\mathrm{spl}}$ denote the midpoint-spline mean error and $\overline E_{\min}$ the smallest mean error among the five practical methods. The percentage excess is
\begin{equation}\label{eq:dos-spline-excess}
    \Delta_{\mathrm{best}}
    :=100\left(\frac{\overline E_{\mathrm{spl}}}{\overline E_{\min}}-1\right).
\end{equation}
Thus $\Delta_{\mathrm{best}}=0$ means that the spline attains the minimum, while a value of $10$ means that its mean error is $10\%$ larger; this is a relative increase in error, not a percentage-point difference. The ideal spline is excluded from the denominator as well as the ranking. These quantities are computed from the reported means in Table~\ref{tab:dos-relative-l2}, not by averaging ranks or error ratios across individual repetitions.

\begin{table}[!htbp]
    \centering
    \small
    \caption{Midpoint-spline rank and percentage excess over the smallest practical-method mean relative $L^2$ error in Table~\ref{tab:dos-relative-l2}. Rank $1$ denotes the smallest mean among the five practical estimators, excluding the ideal-spline reference. Excess values are computed using~\eqref{eq:dos-spline-excess} and rounded to one decimal place. Here $M$ is the number of iterations per probe and $R=5$. These descriptive comparisons do not establish statistical significance.}
    \label{tab:dos-spline-rank-gap}
    \begin{tabular}{@{}lrrrr@{}}
        \toprule
        & \multicolumn{2}{c}{$M=15$} & \multicolumn{2}{c}{$M=90$} \\
        \cmidrule(lr){2-3}\cmidrule(l){4-5}
        Matrix & Spline rank & Excess over best (\%)
               & Spline rank & Excess over best (\%) \\
        \midrule
        \texttt{benzene}       & 1 &  0.0 & 4 & 14.9 \\
        \texttt{HubbardL8}     & 2 &  2.5 & 2 &  7.7 \\
        \texttt{HeisenbergL16} & 3 & 15.9 & 3 &  5.9 \\
        \texttt{Li7Nmax2}      & 3 &  2.6 & 1 &  0.0 \\
        \texttt{nd3k}         & 1 &  0.0 & 2 & 17.4 \\
        \texttt{fv1}          & 1 &  0.0 & 3 &  7.9 \\
        \bottomrule
    \end{tabular}
\end{table}

Across the twelve listed matrix--iteration settings, the midpoint spline ranks among the three most accurate practical estimators in eleven cases; it ranks fourth for \texttt{benzene} at $M=90$. Even in that case, its mean error is only $14.9\%$ above the smallest practical-method mean, and its largest percentage excess over all twelve settings is approximately $17.4\%$. At $M=15$, the three cases in which it does not attain the minimum have excess errors of $2.5\%$, $15.9\%$, and $2.6\%$. At $M=90$, the excess in the five such cases ranges from approximately $5.9\%$ to $17.4\%$. Thus rank alone does not convey the magnitude of a difference. With equal weight given to each setting, the spline has the lowest average rank, $2.17$, compared with $2.83$ for both dual GP and KPM--Jackson, $3.42$ for single GP, and $3.75$ for stochastic Lanczos. The reduction from three first-place results at $M=15$ to one at $M=90$ therefore does not imply that the spline ceases to be competitive.

\paragraph{Comparison with the ideal-spline reference.}
The ideal spline in Table~\ref{tab:dos-relative-l2} serves a different purpose: it indicates the accuracy of a spline construction using spectral information, rather than that of an alternative practical estimator. Its mean error is close to the midpoint spline's at $M=15$ for most matrices. At $M=90$, the gap is more pronounced on some problems; for \texttt{benzene}, the midpoint and ideal spline means are $0.03870$ and $0.01818$, respectively. However, on \texttt{nd3k}, the midpoint spline has a slightly smaller mean at both listed values of $M$. The ideal construction is therefore an informative baseline, not a guaranteed lower error bound. These comparisons alone do not separate the effects of stochastic midpoint error and interpolation.

\paragraph{Practical implications.}
The tables support two practical advantages of the midpoint spline. It gives smaller errors than the tested direct Lanczos and KPM estimators with few iterations per probe, while remaining consistently competitive with the best practical method across both listed iteration counts. Together with its built-in nonnegativity and normalization, this makes the spline a useful choice when the best-performing method is not known in advance. The density plots and other metrics remain essential to this interpretation: mean $L^2$ error does not fully characterize local spectral detail or certify density validity. In particular, GP estimates remain in the rankings even when they contain negative values in some regions of the spectrum, so their norm accuracy must be considered alongside the Jensen--Shannon validity counts.

\paragraph{Additional probe-count tests.}
In additional tests, we held the number of iterations per probe $M$ fixed and varied the number of probes over $R=10,15,20$. These tests showed little improvement in DOS approximation accuracy across the tested methods as $R$ increased. The results are summarized qualitatively here; no additional figures or error tables are included. Thus, for the configurations examined, increasing the probe count alone incurred more matvecs without an appreciable gain in accuracy. This observation does not establish that accuracy is generally insensitive to $R$, or that a small probe count suffices for other matrices and iteration counts.
% Reproducibility note: specify the fixed M value(s) and matrices used in these additional probe-count tests before submission.

%Chao: I would leave the following section out since it may cause confusion and invite criticism. The reviewer may ask for more tests to support the claim, which is not the focus of this paper. We can discuss this in the rebuttal if needed.
%\paragraph{Scope of the evidence.}
%The tabulated comparisons concern the six test matrices at $R=5$ and one selected Gaussian width per matrix; the additional probe-count tests provide a limited qualitative check beyond that setting. Error gaps are evaluated within each matrix at its declared resolution; neither the ranks nor the gaps establish statistical significance or a universal ordering of methods. The experiments do not establish general robustness to probe count, bandwidth, or endpoint choices, asymptotic convergence of the compression induced by averaging corresponding ordered nodes and weights, or an elapsed-time advantage. With dimensions from $1961$ to $12870$, they demonstrate behavior on varied, fully verifiable spectra rather than large-scale scalability. Finally, validation smoothing can attenuate discrepancies in the raw derivatives, so the reported accuracy does not imply recovery of the unsmoothed DOS to arbitrary resolution.

\FloatBarrier

\section{Summary and Conclusions} \label{sec:conclusion}
We have studied DOS estimation by differentiating a smooth approximation to the cumulative spectral distribution. Building on Lanczos quadrature and monotone cumulative interpolation, the proposed midpoint-spline method averages corresponding ordered Ritz values and their weights across runs, fits the resulting jump midpoints, and differentiates the piecewise-cubic fit analytically. The resulting piecewise-quadratic DOS is nonnegative and has unit mass. Its construction requires no Gaussian broadening bandwidth, although the knot distribution and artificial endpoints determine an intrinsic resolution. GP fits provide an alternative with one or two correlation scales and model-conditional derivative uncertainty, but do not by themselves guarantee a normalized, nonnegative density.  It is possible to extend a two-scale GP to a multi-scale GP by adding more correlation kernels, but this increases the number of hyperparameters to be estimated and the computational cost.  The midpoint-spline method is simple to implement, requires no hyperparameter tuning, and has a built-in density constraint.

We compared the midpoint-spline interpolation based approach with the basic Gaussian regularized stochastic Lanczos approach, the KPM--Jackson method and  GP fits on six test problems spanning electronic structure, interacting fermions, disordered spins, nuclear configuration interaction, and molecular vibrational anslysis. Our results support a specific practical advantage: the midpoint spline often recovers useful density information before a Gaussian sum over the available Lanczos nodes is sufficiently resolved. 
%CY: don't need to be so specific in the conclusion
%At $M=15$ and five probes, its mean relative $L^2$ error is lower than those of stochastic Lanczos and KPM--Jackson for every tested matrix. At $M=90$, the benefit remains substantial for the oscillatory Lanczos estimates of \texttt{benzene}, \texttt{HeisenbergL16}, and \texttt{Li7Nmax2}. It is not universal: stochastic Lanczos, KPM--Jackson, or a GP fit can attain smaller errors on other problems or under other metrics. The ideal spline serves as an informative reference, not a guaranteed lower error bound.

The complementary error measures distinguish peak-height errors, integrated discrepancies, and normalized shape differences. In particular, good norm accuracy of a signed GP estimate does not make it a valid density, and undefined Jensen--Shannon values must be reported alongside the valid-repetition statistics. All numerical comparisons use a common Gaussian validation filter; the conclusions concern that declared resolution, not pointwise recovery of the exact distribution-valued DOS. Likewise, equal nominal counts of matrix--vector products do not establish equal elapsed time.

The main unresolved questions concern the accuracy and stability of the compression induced by averaging corresponding ordered nodes and weights, the effect of endpoint extrapolation, and the propagation of cumulative-fit error through differentiation. Further work should extend the additional probe-count tests to a systematic sensitivity study, examine the effect of the validation bandwidth, compare compression with fitting the full pooled quadrature data, and assess time and memory costs on larger problems. Constrained GP models and calibration of their uncertainty estimates provide additional directions. Within the tested setting, the midpoint spline offers a simple combination of accuracy with few matvecs and built-in density constraints, making it a useful alternative when the objective is an informative spectral profile rather than individual eigenvalue resolution.

\section*{Acknowledgement}
This work was supported in part by the U.S. Department of Energy, Office of Science, Office of Advanced Scientific Computing Research, Scientific Discovery through Advanced Computing (SciDAC) Program through the FASTMath Institute and in part by the U.S. Department of Energy, Office of Science, Office of Advanced Scientific Computing Research's Applied Mathematics Competitive Portfolios program under U.S. Department of Energy Contract No. DE-AC02-05CH11231. This work used computing resources of the National Energy Research Scientific Computing Center (NERSC), which is supported by the DOE Office of Science under Contract No. DE-AC02-05CH11231, using NERSC Awards ASCR-ERCAP m1027 for 2026.

\bibliographystyle{plain}
\bibliography{reference}

\end{document}